\documentclass[12pt]{article}
\usepackage[utf8]{inputenc}
\usepackage[T1]{fontenc}
\usepackage{lmodern}
\usepackage[margin=2.5cm]{geometry}
\usepackage[english]{babel}
\usepackage{amsmath, amssymb, amsthm}
\usepackage{amsrefs}
\usepackage{authblk}
\usepackage{hyperref}
\usepackage[all]{xy}
\usepackage{tikz}
\usepackage{xypic}
\usepackage{xcolor}

\newtheorem{thm}{Theorem}[section]
\newcounter{thmM}
\newtheorem*{thml}{Theorem \thethmM \stepcounter{thmM}}

\newtheorem{conj}[thm]{Conjecture}

\newtheorem{cor}[thm]{Corollary}

\newtheorem{defi}[thm]{Definition}

\newtheorem{lem}[thm]{Lemma}

\newtheorem{prop}[thm]{Proposition}

\newenvironment{demo}{\noindent \textit{Proof.}}{\hfill\QED\\}

\newenvironment{rmq}{\stepcounter{thm}\noindent \textbf{Remark \thethm}}{}
\newenvironment{rmqs}{\stepcounter{thm}\noindent \textbf{Remarks \thethm}}{}

\renewcommand*{\bar}[1]{\overline{\vphantom{#1^I}#1}}
\renewcommand*{\Re}{\mathfrak{Re}}

\newcommand*{\C}{\mathbb C}
\newcommand*{\E}{\mathbb E}

\renewcommand*{\H}{\mathbb H}
\newcommand*{\K}{\mathbb K}

\renewcommand*{\O}{\mathcal O}
\renewcommand*{\P}{\mathbb P}
\newcommand*{\Q}{\mathbb Q}
\newcommand*{\R}{\mathbb R}
\renewcommand*{\S}{\mathbb S}

\newcommand*{\Z}{\mathbb Z}
\newcommand*{\I}{\mathcal I}

\newcommand*{\e}{\mathrm e}
\newcommand*{\z}{\mathfrak o}
\renewcommand*{\I}{\mathcal I}

\newcommand*{\QED}{\hfill\ensuremath{\square}}
\newcommand*{\trig}{\, \triangleleft \,}
\newcommand*{\p}{\mathfrak p}
\newcommand*{\Pp}{\mathfrak P}

\DeclareMathOperator{\Var}{Var}
\DeclareMathOperator{\ord}{ord}
\DeclareMathOperator{\Li}{Li}
\DeclareMathOperator{\Irr}{Irr}
\DeclareMathOperator{\Irrr}{Irr^{real}}

\DeclareMathOperator{\Ind}{Ind}
\DeclareMathOperator{\reg}{reg}
\DeclareMathOperator{\Gal}{Gal}
\DeclareMathOperator{\codim}{codim}

\DeclareMathOperator{\GL}{\bf GL}

\DeclareMathOperator{\BM}{\mathsf{BM}}
\DeclareMathOperator{\GRH}{\mathsf{GRH}}
\DeclareMathOperator{\LI}{\mathsf{LI}}
\DeclareMathOperator{\LIP}{\mathsf{LI^+}}
\DeclareMathOperator{\LIM}{\mathsf{LI^-}}

\begin{document}

\title{\textbf{Chebyshev's bias in dihedral and generalized quaternion Galois groups}}
\author{Alexandre Bailleul}
\date{}

\maketitle

\begin{abstract} We study the inequities in the distribution of Frobenius elements in Galois extensions of the rational numbers with Galois groups that are either dihedral $D_{2^{n}}$ or (generalized) quaternion $\H_{2^n}$ of two-power order. In the spirit of recent work of Fiorilli and Jouve \cite{FiJo}, we study, under natural hypotheses, some families of such extensions, in a horizontal aspect, where the degree is fixed, and in a vertical aspect, where the degree goes to infinity. Our main contribution uncovers in families of extensions a phenomenon, for which Ng gave numerical evidence in \cite{Ng} : real zeros of Artin $L$-functions sometimes have a radical influence on the distribution of Frobenius elements.
\end{abstract}

\subsection*{Introduction}

The prime number theorem in arithmetic progressions states that if $q \geq 1$ is an integer and $a$ is prime to $q$, then $$\pi(x,q,a) := \left|\{p \leq x \mid p = a \text{ mod } q\}\right| \underset{x \to +\infty}{\sim} \frac{1}{\varphi(q)} \Li(x),$$ where $\Li$ is the logarithmic integral. Even though for $a$ and $b$ coprime to $q$, $\pi(x, q, a)$ and $\pi(x,q,b)$ have the same asymptotic value, it may happen that one count could be larger than the other most of the time. This is Chebyshev's bias as observed initially in $1853$ in \cite{Cheb} : there seems to be more primes congruent to $3$ mod $4$ than $1$ mod $4$ in generic intervals $[2, x]$.

These so-called "prime number races" have been studied extensively by many authors, most notably by Knapowski and Tur\'{a}n in two series of papers \cite{KT1} and \cite{KT2}, and more recently by Kaczorowski (\cite{Kac}, \cite{Kac2}) and Rubinstein-Sarnak in \cite{RuSa}. The latter managed to explain conditionally Chebyshev's bias : under natural hypotheses which we will consider later, the inequality $\pi(x,4,3) > \pi(x,4, 1)$ holds "$99.59\%$ of the time" (this will be explained rigorously in a later paragraph). For a general modulus $q$, they have shown a bias towards non-quadratic residues mod $q$ against quadratic residues mod $q$. This lead to a lot of subsequent works which either aimed at giving asymptotic formulas for the proportion of time during which one of the prime number team is in the lead (\cite{FiMa}), weakening some of the aforementionned hypotheses of Rubinstein and Sarnak (\cite{MN}), studying races with many participants (\cite{Lam1}, \cite{Lam2}, \cite{Lam3}, \cite{Lam4}), or generalising Rubinstein and Sarnak's work to other contexts (\cite{Cha}, \cite{ChaIm}, \cite{DevMen}, \cite[Chapter 5]{Ng}).

In \cite{Ng}, Ng, following a suggestion made in \cite{RuSa}, extended Rubinstein and Sarnak's framework to conjugacy classes of the Galois group $G$ of a Galois extension $L/K$ of number fields, in the context of the Chebotarev density theorem (\cite[§2, Théorème 1]{SerCheb}). Recall that if $C$ is such a conjugacy class, then $$\pi(x, C, L/K) := \left|\left\{\p \trig \mathcal{O}_K \text{ unramified} \mid N(\p) \leq x, \left(\frac{\p}{L/K}\right) = C\right\}\right| \underset{x \to +\infty}{\sim} \frac{|C|}{|G|} \Li(x),$$ where $N(\p)$ is the norm of the prime ideal $\p$ and $\left(\frac{\p}{L/K}\right)$ is the Artin symbol, or Frobenius conjugacy class, at $\p$. Taking $L=\Q(\zeta_q)$, where $\zeta_q$ is a primitive $q$-th root of unity in $\C$, and $K=\Q$, we get the prime number theorem in arithmetic progressions and one recovers the setting of races between prime numbers in arithmetic progressions of \cite{RuSa}, as it is immediate to see that the Frobenius automorphism associated with the prime $p$ only depends on the congruence class of $p$ modulo $q$. We say we study \emph{Chebotarev biases} when comparing the behaviours of prime ideal counting functions $\pi(x, C_1, L/K)$ and $\pi(x, C_2, L/K)$ for distinct conjugacy classes $C_1, C_2$ of $\Gal(L/K)$.

In recent work, Fiorilli and Jouve (\cite{FiJo}) have shown some Chebotarev races to be extremely biased (\textit{i.e.} as for the original case where one compares $\pi(x, 4, 3)$ and $\pi(x, 4, 1)$, the underlying density is close to $1$) by using large deviations principles. In the opposite direction, they managed to show central limit theorem behaviours which correspond to moderately biased Chebotarev races (\textit{i.e.} the underlying density is close to $\frac{1}{2}$). In both cases, the asymptotic results appear as we let the conductors go to infinity. In \cite{FiJo}, the theoretical results are applied to some families of number field extensions with Galois groups dihedral of order $2p$ ($p$ an odd prime), quasidihedral, symmetric or Frobenius of order $p(p-1)$.

In this paper, we exploit work of Fröhlich (\cite{Fro1}, \cite{Fro2}, \cite{Fro3}) on root numbers of quaternion extensions and we perform class field theoretic constructions to study, in the context of \cite{FiJo}, Chebyshev's bias in some families of number field extensions with Galois groups  of $2$-power order : dihedral or (generalized) quaternion. Our main contributions highlight the role of central zeros of Artin $L$-functions in the study of this bias. They come as the result of the two following points of view. We first work in a "horizontal aspect", in which the Galois group is fixed (up to isomorphism). In this context, thanks to a result of Fröhlich (Theorem \ref{Fro1}), we construct two different kinds of families of number fields which have quaternion Galois groups over $\Q$, but in which the existence or not of a central zero for the corresponding Dedekind zeta function leads to opposite biases (Theorem \ref{Opp8}). Second, in the "vertical aspect", we build arbitrarily high degree towers of number field extensions with Galois groups dihedral or quaternion of $2$-power order, and exhibit different kinds of behaviours in both families (Theorems \ref{TabD} and \ref{TabQ}). We obtain extreme bias, moderate bias and unbiased races (the underlying density is plainly equal to $\frac{1}{2}$). We emphasize that some of these behaviours are directly linked to the existence or not of a central zero for the corresponding Dedekind zeta functions in the quaternion case. This provides the first theoretical construction that confirms numerical evidence obtained by Ng. Finally, we are able to prove partial monotonicity phenomena in the evolution of the bias inside the tower itself (Theorem \ref{Mono}).

\subsubsection*{Outline of the paper}

In Section \ref{Sect1} we introduce the preliminary results and conjectures needed to study Chebyshev's bias in number field extensions, and state abridged versions of our main results. In Section \ref{Conduc} we recall the definition of the Artin conductor of a complex character of the Galois group of a number field extension, and relate it to the study of Chebyshev's bias in the extension considered. Then in Section \ref{Sect2.2}, we study the character-theoretic properties of dihedral and generalized quaternion groups of $2$-power order. In Section \ref{Sect3} we focus on the "horizontal aspect" of our study, in which the Galois group is fixed to be a quaternion group of order $8$. In Sections \ref{Sect4.1} and \ref{Sect4.2} we give bounds on moments of the random variables governing Chebyshev's bias in Galois extensions with Galois groups dihedral and generalized quaternion groups of $2$-power order. In Section \ref{Sect4.3} we construct such extensions with controlled ramification, and in Section \ref{Sect4.4} we give estimates on Chebyshev's bias in those extensions.

\section{Chebotarev biases in Galois groups of number fields}

\label{Sect1}

\subsection{Notations and main results}

We recall the following standard notations.\\

In this section, we let $L/K$ be a Galois number field extension with Galois group $G$. Let $C_1$ and $C_2$ be distinct conjugacy classes of $G$. We would like to give a notion of measure to the set $$\mathcal{P}_{L/K, C_1, C_2} := \left\{x \geq 2 \mid \frac{\pi(x, C_1, L/K)}{|C_1|} > \frac{\pi(x, C_2, L/K)}{|C_2|}\right\},$$ where $C_1$ and $C_2$ are conjugacy classes of $G$. A first guess would be to use the natural density, defined as the limit as $x$ tends to $+\infty$ of $$\frac{\left|\mathcal{P}_{L/K, C_1, C_2} \cap [0, x]\right|}{x},$$ where $|.|$ denotes the Lebesgue measure. Unfortunately, it has been shown that this limit does not exist, even in the case of prime number races \cite{Kac} (though one should note that for races for irreducible polynomials over finite fields studied in \cite{Cha} for instance, the corresponding natural density exists). Therefore, we use the notion of logarithmic density, better-suited to the study of Chebyshev's bias as was first observed in \cite{Win}.

\begin{defi} Let $A$ be a Borel set in $\R$. Its logarithmic density is, when it exists, $$\delta(A) := \lim_{X \to +\infty} \frac{1}{X} \int_2^X \mathbf{1}_A(\e^t) \,dt.$$ If $C_1$ and $C_2$ are conjugacy classes of $G$, we write $\delta(L/K, C_1, C_2)$ for $\delta(\mathcal{P}_{L/K, C_1, C_2})$. We say the race between $C_1$ and $C_2$ is unbiased if $\delta(L/K, C_1, C_2) = \frac{1}{2}$, that it is biased towards $C_1$ if $\delta(L/K, C_1, C_2) > \frac{1}{2}$ and biased towards $C_2$ if $\delta(L/K, C_1, C_2) < \frac{1}{2}$.
\end{defi}

Ng has shown in \cite{Ng} that this logarithmic density always exists under some natural hypotheses (see Section \ref{Conj}). Our goal is to give estimates for such densities.\\

We now state our main results. The hypothesis $\GRH, \LI^-, \LI$ and $\LI^+$ will be explained in details in Section \ref{Conj}. Our first result is in a horizontal aspect : we exhibit moderately biased races in families of extensions of $\Q$ with fixed Galois group isomorphic to $\H_8$, the usual quaternion group of order $8$.

\begin{thml}\hypertarget{ThA} Assume $\GRH$ and $\LIP$. For any function $f$ such that $f(n) \underset{n \to +\infty}{\longrightarrow} +\infty$, there exist two families $(K_d)_d$ and $(L_d)_d$ of number fields, indexed by square-free integers satisfying $d > 1$ and $d=1$ mod $4$, such that for any $d$ in the index set :

\begin{itemize}
\item[i)] $\Q(\sqrt d) \subset K_d \cap L_d$.
\item[ii)] $\Gal(K_d/\Q) \simeq \Gal(L_d/\Q) \simeq \H_8$.
\item[iii)] $0 < \frac{1}{2} - \delta(K_d/\Q, C_1, C_{-1}) \ll \frac{1}{f(d)}$, where $C_1$ and $C_{-1}$ denote the conjugacy classes of $1$ and $-1$ in $\H_8$.
\item[iv)] $0 < \delta(L_d/\Q, C_1, C_{-1}) - \frac{1}{2} \ll \frac{1}{f(d)}$.
\end{itemize}
\end{thml}

The most important feature of this result is that we are able to exhibit two families in which the biases are opposite to each other, one is $< \frac{1}{2}$ while the other is $> \frac{1}{2}$, the difference coming from the existence of a central zero for the corresponding Dedekind zeta function in the first case, and the absence of such a zero in the second case. Moreover, the prime number races considered can be specified to be as moderately biased as possible, in the sense that the logarithmic densities can be made arbitrarily close to $\frac{1}{2}$ when $d$ grows, while ensuring that $\Q(\sqrt{d}) \subset K_d \cap L_d$. In fact the choice of $f$ allows one to have arbitrarily large variances for the random variables $X(K_d/\Q, C_1, C_{-1})$ and $X(L_d/\Q, C_1, C_{-1})$ governing the biases (see Theorem \ref{VA}).\\

The second part of our work is in a vertical aspect : we build families of extensions of $\Q$ with Galois groups dihedral and quaternion of $2$-power order with arbitrarily large degree, in which we are able to show extreme biases, moderate biases and absence of bias.

\begin{thml}\hypertarget{ThB} Assume $\GRH$ and $\LIP$. There exist absolute constants $c_1, c_2, c_3 > 0$ and a family $(\mathcal{D}_n)_n$ of number fields such that for any $n \geq 3$, the following hold :

\begin{itemize}
\item[i)] $\Gal(\mathcal{D}_n/\Q) \simeq D_{2^{n-1}} = \langle r, s \mid r^{2^{n-1}}=s^2=1, srs^{-1}=r^{-1} \rangle$.
\item[ii)] $c_1\exp(-c_22^n) < \delta(\mathcal{D}_n/\Q, C_1, C_{-1}) < \exp\left(-c_3 \frac{2^n}{n}\right)$ where $C_a$ denotes the conjugacy class of $a$, and $-1$ denotes $r^{2^{n-2}} \in D_{2^{n-1}}$.
\item[iii)] $0 < \frac{1}{2} - \delta(\mathcal{D}_n/\Q, C_{-1}, C_s) \ll \frac{1}{2^{n/3}}$.
\item[iv)] $\delta(\mathcal{D}_n/\Q, C_{r^k}, C_{s}) = \delta(\mathcal{D}_n/\Q, C_{r^k}, C_{rs}) = \frac{1}{2}$ when $1 \leq k \leq 2^{n-2}-1$ is odd.
\end{itemize}
\end{thml}

We chose a sample of possible couples of conjugacy classes in the above statement. For an exhaustive treatment of conjugacy classes, see Theorem \ref{TabD}. In \cite{FiJo}, Fiorilli and Jouve, relying on a construction of Klüners, only deal with dihedral groups of order $2p$ with $p$ varying in the set of odd primes. The following result deals with generalized quaternion Galois groups.

\begin{thml}\hypertarget{ThC} Assume $\GRH$ and $\LIP$. There exist absolute constants $c_1, c_2, c_3 > 0$ and two families $(\mathcal Q_n^+)_n$ and $(\mathcal{Q}_n^-)_n$ of number fields such that for any $n \geq 3$, the following hold :

\begin{itemize}
\item[i)] $\Gal(\mathcal Q_n^{\pm}/\Q) \simeq \H_{2^n} = \langle x, y \mid x^{2^{n-1}}=1, x^{2^{n-2}}=y^2, yxy^{-1}=x^{-1} \rangle$, and the root number of each symplectic character of $\Gal(\mathcal Q_n^{\pm}/\Q)$ is the same (we denote it by $W_{\mathcal Q_n^{\pm}}$). Moreover, $W_{\mathcal{Q}_n^+}=1$ and $W_{\mathcal{Q}_n^-} = -1$.
\item[ii)] $c_1\exp(-c_22^n) < \left|\frac{1-W_{\mathcal Q_n^{\pm}}}{2} - \delta(\mathcal Q_n^{\pm}/\Q, C_1, C_{-1})\right| < \exp\left(-c_3 \frac{2^n}{n}\right),$ where $C_a$ denotes the conjugacy class of $a$, and $-1$ denotes $x^{2^{n-2}}$.
\item[iii)] $\delta(\mathcal Q_n^+/\Q, C_{1}, C_{x^{k}}) = \frac{1}{2}$ when $1 \leq k \leq 2^{n-2}-1$ is even and $c_1\exp(-c_2 32^n) < \delta(\mathcal Q_n^-/\Q, C_{1}, C_{x^k}) < \exp\left(-c_3 \frac{2^n}{n}\right)$ for $1 \leq k \leq 2^{n-2}-1$.
\item[iv)] $0 < \frac{1}{2} - \delta(\mathcal Q_n^+/\Q, C_1, C_y) \ll \frac{1}{2^{n/3}}$ and $c_1\exp(-c_2 2^n) < \delta(\mathcal Q_n^-/\Q, C_{1}, C_{y}) < \exp\left(-c_3 \frac{2^n}{n}\right)$.
\end{itemize}
\end{thml}

Again, we chose to focus on a sample of meaningful cases of bias estimates in ii), iii) and iv), but more have been computed and are stated in Theorem \ref{TabQ}. The new and remarkable feature of this theorem is that it highlights the role played by the existence or not of a central zero for the corresponding Dedekind zeta function. The presence (or the inexistence) of such a zero leads to : extreme biases of opposite signs in ii), extreme biases or no biases in iii), extreme bias or moderate bias in iv).\\

Finally, our last contribution is to show that we can observe a monotonicity phenomenon in the evolution of Chebyshev's bias in subextensions of dihedral extensions of $\Q$.

\begin{thml}\hypertarget{ThD} Assume $\GRH$ and $\LI$. Consider the number fields $(\mathcal{D}_n)_n$ as in Theorem \hyperlink{ThB}{B}. For every $n \geq 3$, there are dihedral subfields $\Q = \mathcal{D}_n^{(n)} \subset \mathcal{D}_n^{(n-1)} \subset \dots \subset \mathcal{D}_n^{(3)} \subset \mathcal{D}_n$ such that $\Gal(\mathcal{D}_n/\mathcal{D}_n^{(i)}) \simeq D_{2^n}$ for $3 \leq i \leq n$, and such that for any $\varepsilon > 0$ and any sufficiently large $n$, for $3 \leq i < j \leq n$ satisfying $i \leq n\frac{1 + \varepsilon}{2}$ and $j \geq n\left(\frac{1 + 3\varepsilon}{2}\right)$, we have $$\delta(\mathcal{D}_n/\mathcal{D}_n^{(j)}, C_1^{(j)}, C_{-1}^{(j)}) < \delta(\mathcal{D}_n/\mathcal{D}_n^{(i)}, C_1^{(i)}, C_{-1}^{(i)}),$$ where $C_1^{(k)}$ and $C_{-1}^{(k)}$ denote the conjugacy classes of $1$ and $-1$ in $\Gal(\mathcal{D}_n/\mathcal{D}_n^{(k)})$.
\end{thml}

This statement means that, in the subextensions of $\mathcal{D}_n/\mathcal{D}_n^{(i)}$, Chebyshev's bias is more extreme at the bottom of the tower than at the top. A more general result including the quaternion fields $\mathcal Q_n^{\pm}$ of Theorem \hyperlink{ThC}{C} is stated as Theorem \ref{Mono}. The proof relies on a new large deviation bound for the values of Chebyshev's bias in families of number field extensions (see Theorem \ref{Dev}).

\subsection{Recollection on Artin characters}

A certain class of characters of $G$ plays a particular role in our study, the class of symplectic characters. To introduce them, we need the notion of Frobenius-Schur index of a character.

\begin{defi} Let $\chi$ be a character of $G$. We define its Frobenius-Schur index by $$\varepsilon_2(\chi) := \frac{1}{|G|} \sum_{g \in G} \chi(g^2).$$ We also define $\Irr(G)$ to be the set of irreducible characters of $G$, and $\Irrr(G)$ to be the set of real-valued irreducible characters of $G$.
\end{defi}

The Frobenius-Schur index of an irreducible (complex) character of $G$ determines if such a character can be afforded by a real-valued representation or not, thanks to the following well-known result from character theory.

\begin{thm}[\cite{Isa} p.58]\label{Frob} Let $\chi$ be an irreducible complex character of $G$. Then only one of the following three statements holds :

\begin{enumerate}
\item[i)] $\varepsilon_2(\chi)=0$, in this case $\chi$ is not real-valued. We say that $\chi$ is unitary.
\item[ii)] $\varepsilon_2(\chi)=1$, in this case $\chi$ is real-valued and can be afforded by a real-valued representation of $G$. We say that $\chi$ is orthogonal.
\item[iii)] $\varepsilon_2(\chi)=-1$, in this case $\chi$ is real-valued and cannot be afforded by a real-valued representation of $G$. We say that $\chi$ is symplectic.
\end{enumerate}
\end{thm}

It is expected that symplectic characters are exactly the irreducible characters which can yield real zeros for their associated Artin $L$-function (see conjecture $\LI$ below).

Finally, we recall that if $\chi$ is an irreducible character of $G$, then its Artin $L$-function satisfies a functional equation (\cite[p.28]{MuMu}) of the form $$\Lambda(s, \chi, L/K) = W(\chi, L/K) \Lambda(1-s, \overline{\chi}, L/K),$$ where $W(\chi, L/K)$ is a complex number of modulus $1$, called the root number of $\chi$, and $s \mapsto \Lambda(s, \chi, L/K)$ is the completed Artin $L$-function associated to $\chi$, which is the product of $s \mapsto L(s, \chi, L/K)$ with Gamma factors coming from archimedian places. Unless there is an ambiguity in the extension considered, we will usually write $W(\chi)$ for $W(\chi, L/K)$. The root number satisfies $\bar{W(\chi)} = W(\overline{\chi})$. In particular if $\chi$ is real-valued then $W(\chi) = \pm 1$. We have the following important way to detect central zeros of Artin $L$-functions.

\begin{prop}\label{Zerr} If $\chi \in \Irrr(G)$ with $W(\chi)=-1$, then $L\left(\frac{1}{2}, \chi, L/K\right) = 0$.
\end{prop}

\begin{demo} Since $\chi$ is real-valued, evaluating the functional equation at $\frac{1}{2}$ we find $$\Lambda\left(\frac{1}{2}, \chi \right) = - \Lambda\left(\frac{1}{2}, \chi \right),$$ i.e. $\Lambda\left(\frac{1}{2}, \chi \right) = 0$. Since the Gamma function never vanishes on $\C$, this implies that $L\left(\frac{1}{2}, \chi, L/K \right)=0$.
\end{demo}

It is expected that the converse also holds for real-valued characters when the base field is $\Q$, see conjecture $\LIP$ below.

\subsection{Conjectures}

\label{Conj}

As in \cite{FiJo}, we consider natural conjectures on the distribution of zeros and poles of Artin $L$-functions, some of which are generalizations of conjectures used in the work of Rubinstein and Sarnak. Recall that $L/K$ is a Galois extension of number fields.

\begin{conj}[Artin's conjecture] If $\chi$ is a non-trivial irreducible character of $G$, then $s \mapsto L(s, \chi, L/K)$ is entire.
\end{conj}

\begin{conj}[$\GRH$] If $\chi$ is an irreducible character of $G$, then the non-trivial zeros of $s \mapsto L(s, \chi, L/K)$ have real part $\frac{1}{2}$.
\end{conj}

\begin{conj}[$\LIM$] Let $L_0/\Q$ be the Galois closure of $L/\Q$. Then the multiset of imaginary parts of zeros $$\Gamma_{L_0/\Q} := \bigcup_{\chi \in \Irr(\Gal(L_0/\Q))} \left\{\gamma > 0 \mid  L\left(\frac{1}{2} + i \gamma, \chi, L_0/\Q\right)=0\right\}$$ is linearly independent over $\Q$.
\end{conj}

\begin{conj}[$\LI$] $\LIM$ is true and if $\chi \neq \chi_0$ is a unitary or orthogonal character of $\Gal(L_0/\Q)$ (see Theorem \ref{Frob}) then $L(\frac{1}{2}, \chi, L_0/\Q) \neq 0$. If $\chi$ is a symplectic character of $\Gal(L_0/\Q)$ then $\ord_{s=1/2} L(s, \chi, L_0/\Q)$ is bounded by some absolute constant $M_0$.
\end{conj}

\begin{conj}[$\LIP$] $\LI$ is true and if $\chi$ is a symplectic irreducible character of $\Gal(L_0/\Q)$ then $\ord_{s=1/2} L(s, \chi, L_0/\Q) = \frac{1-W(\chi)}{2}.$
\end{conj}

Let us make a few comments about these conjectures.

\begin{itemize}
\item[•] First, Artin's conjecture needs to be assumed in order to be able to prove explicit formulas for number field analogs of Chebyshev's $\psi$ function, involving sums over zeros of Artin $L$-functions, and discarding the existence of possible poles. In the cases we will consider, namely when $G$ is a dihedral group or generalized quaternion group, Artin's conjecture is known to be true, because such groups are supersolvable, which means they have a normal series $$\{1\} = G_1 \trig G_2 \trig \dots \trig G_k = G$$ where each quotient $G_{i+1}/G_i$ is cyclic for $1 \leq i \leq k-1$. This implies that their irreducible characters are induced by those of abelian subgroups. By inductive properties of Artin $L$-functions, we are reduced to knowing Artin's conjecture in the case of characters of abelian Galois groups, but this is exactly one of the consequences of class field theory, together with work of Hecke on the $L$-functions bearing his name.

\item[•] The Generalized Riemann Hypothesis ($\GRH$) needs to be assumed so that, for conjugacy classes $C_1$ and $C_2$ of $G$, $$\frac{|G|}{|C_2|}\pi(x, C_1, L/K) - \frac{|G|}{|C_2|}\pi(x, C_2, L/K)$$ is oscillating with amplitude of size roughly $\sqrt{x}$. This is central in the analysis of Chebyshev's bias. One can still show the existence of a limiting logarithmic distribution for a convenient renormalization of $\frac{|G|}{|C_2|}\pi(x, C_1, L/K) - \frac{|G|}{|C_2|}\pi(x, C_2, L/K)$ without assuming GRH (see \cite{Dev}), though it depends on the supremum of real parts of non-trivial zeros of the Artin $L$-functions considered.

\item[•] The hypothesis $\LIM$ actually contains two statements. The most obvious one is the linear independence of the positive imaginary parts of zeros of Artin $L$-functions. This hypothesis appears because, in order to understand Chebyshev's bias using explicit formulas for prime counting functions, one needs good information on the joint distribution of the values of $\e^{i \gamma y_1}, e^{i \gamma_2 y}, \dots$ in $\S^1$, where $\gamma_1, \gamma_2, \dots$ denotes the aforementioned positive imaginary parts of zeros, and $y$ varies in $\R$. Again, weaker hypotheses have been used to show the existence of a limiting logarithmic distribution (\cite{Dev}) and results on the bias. Another aspect of $\LI$ we highlight is that it is stated for $L$-functions over $\Q$ instead of $L$-functions over $K$. Indeed, those $L$-functions over $K$ factorize as products of $L$-functions relative to $L_0/\Q$, so linear independence is typically false. The last aspect adressed by $\LI$ is about the multiplicity of the zeros of Artin $L$-functions over $\Q$. One should expect those Artin $L$-functions to be "primitive" in the sense of \cite{RudSa}, and those should not satisfy any kind of non-trivial algebraic relations, except for their respective functional equations. Instead of the simplicity of the zeros, we could have opted for an hypothesis about the boundedness of the multiplicities of such zeros (called $\BM$ in \cite{FiJo}), which would have been enough for a few intermediate results of this paper.

\item[•] The assumptions on the order of vanishing at $\frac{1}{2}$ in $\LI$ and $\LIP$ appear because the quantity $z(C)$ defined above is involved in the mean of the random variables attached to the limiting distribution governing the bias. Examples of symplectic characters with root numbers $-1$ were first given by Armitage (\cite{Arm}) and Serre (unpublished). In particular, they yield Artin $L$-functions vasnishing at $1/2$ because of Proposition \ref{Zerr}. As was pointed out by the author of the present paper to the authors of \cite{FiJo} while both papers were under preliminary form, using one of the aforementionned examples, the conjecture "Artin $L$-functions attached to orthogonal or unitary characters do not vanish at $1/2$" could not hold over general number fields $K$. Indeed, consider Serre's example, as described in \cite[section 5.3.3]{Ng} : let $L=\Q(\theta)$, where $\theta = \sqrt{\frac{5+\sqrt{5}}{2} \frac{41 + \sqrt{5 \cdot 41}}{2}}$. It is a Galois extension of $\Q$ with Galois group isomorphic to the quaternion group $\H_8$ (see section 2). The root number of its non-abelian (symplectic) character $\psi$ is $-1$, so $L(1/2, \psi, L/\Q)=0$. Now consider the subfield $K=\Q(\sqrt 5) \subset L$. We have $\Gal(L/K) \simeq \Z/2\Z \times \Z/2\Z$ which does not admit any irreducible symplectic character. But we have the classical factorization \begin{align*}
\zeta_L(s) &= P_1(s) L(s, \psi, L/\Q)^2\\ &= P_2(s),
\end{align*} where $P_1$ is the product of the Artin $L$-functions attached to the four abelian irreducible characters of $\Gal(L/\Q)$, and $P_2$ is the product of the Artin $L$-functions attached to the four irreducible characters of $\Gal(L/K)$. Since $L(1/2, \psi, L/\Q)=0$, at least one of those Artin $L$-functions, attached to a non-symplectic irreducible character, must vanish at $1/2$. This shows that we cannot expect easy non-vanishing statements at the central point in the relative case $L/K$ when $K \neq \Q$, and that such a statement should involve the way irreducible characters of $G$ are "induced" to Artin characters over $\Q$.
\end{itemize}

\subsection{The probabilistic approach to Chebotarev biases}

Let us recall the setting of \cite{FiJo} together with useful results. As before, $L/K$ is a Galois extension of number fields with group $G$. We will assume $L/\Q$ is Galois, with Galois group $G^+$, so that $G$ is a subgroup of $G^+$. This will be the case in our applications. If $C$ is a conjugacy class of $G$, then we denote by $C^+$ the conjugacy class generated by $C$ in $G^+$, \textit{i.e.} $C^+ := \bigcup_{a \in G^+} aCa^{-1}$. In what follows, if $F/E$ is a Galois extension of number fields and $\chi \in \Irr(\Gal(F/E))$, a summation over $\gamma_{\chi} > 0$ means a summation over the corresponding zero multiset $$\Gamma_{F/E, \chi} := \{\gamma > 0 \mid  L(1/2 + i \gamma, \chi, F/E)=0\}.$$ We also write $\Gamma_{F/E} := \bigcup_{\chi \in \Irr(\Gal(F/E))} \Gamma_{F/E, \chi}$ and $\Gamma^{\text{real}}_{F/E} := \bigcup_{\chi \in \Irrr(\Gal(F/E))} \Gamma_{F/E, \chi}$.\\

Following \cite{FiJo}, we introduce some quantities related to the conjugacy classes under consideration, but also to the irreducible characters of $G$.

\begin{defi} Let $C$ be a conjugacy class of $G$. Then we define $$C^{1/2} := \{g \in G \mid g^2 \in C\}$$ and $$z(C) := 2\sum_{\chi \neq \chi_0} \chi(C) \ord_{s=1/2} L(s, \chi, L/K),$$ where the sum is taken over all the non-trivial irreducible characters $\chi$ of $G$, and $s \mapsto L(s, \chi, L/K)$ is the Artin $L$-function associated to $\chi$.
\end{defi}

We summarize in the next statement the main result giving a probabilistic interpretation of the logarithmic densities we are studying.

\begin{thm}[\cite{FiJo}, Proposition 3.18 and Lemma 3.20]\label{VA} Assume Artin's conjecture, $\GRH$ and $\LIM$. For any $\gamma \in \Gamma_{L/\Q}$, we introduce the random variable $X_{\gamma} = \Re(Z_{\gamma})$, where $(Z_{\gamma})_{\gamma}$ is a family of independent random variables uniform on the unit circle. Then for any conjugacy classes $C_1$ and $C_2$ of $G$, we have $$\delta(L/K, C_1, C_2) = \P(X(L/K, C_1, C_2) > 0)$$ where $$X(L/K, C_1, C_2) := \frac{|C_2^{1/2}|}{|C_2|} - \frac{|C_1^{1/2}|}{|C_1|} + z(C_2) - z(C_1) + 2 \sum_{\lambda \in \Irr(G^+)} |\lambda(C_2^+)-\lambda(C_1^+)| \sum_{\gamma_{\lambda} > 0} \frac{X_{\gamma_{\lambda}}}{\sqrt{\frac{1}{4} + \gamma_{\lambda}^2}},$$ unless $C_1^+=C_2^+$.

Moreover, under the previous condition, we have $$\E(X(L/K, C_1, C_2)) = \frac{|C_2^{1/2}|}{|C_2|} - \frac{|C_1^{1/2}|}{|C_1|} + z(C_2) - z(C_1)$$ and 

$$\Var(X(L/K, C_1, C_2)) = 2 \sum_{\lambda \in \Irr(G^+)} |\lambda(C_1^+)-\lambda(C_2^+)|^2 \sum_{\gamma_{\lambda} > 0} \frac{1}{\frac{1}{4} + \gamma_{\lambda}^2}.$$
\end{thm}

Theorem \ref{VA} follows from \cite[Proposition 3.18]{FiJo} and \cite[3.20]{FiJo} by choosing the the class function $t = \frac{|G|}{|C_1|} \mathbf 1_{C_1} - \frac{|G|}{|C_2|} \mathbf 1_{C_2}$, which satisfies $$\hat{t^+}(\lambda) = \lambda(C_1) - \lambda(C_2)$$ for any $\lambda \in \Irr(G^+)$, $$\ord_{s=1/2}L(s, L/K, t) = z(C_1) - z(C_2)$$ and $$\langle t, r \rangle_G = \frac{|C_1^{1/2}|}{|C_1|} - \frac{|C_2^{1/2}|}{|C_2|}.$$

Note that the expression of $\Var(X(L/K, C_1, C_2))$ involves zeros of $L$-functions attached to characters of $G^+$, not of $G$. This is because our linear independence and simplicity hypothesis $\LIM$ was stated over $\Q$. If $C_1^+=C_2^+$ then it is easy to see that $\frac{|C_2^{1/2}|}{|C_2|} - \frac{|C_1^{1/2}|}{|C_1|} + z(C_2) - z(C_1)=0$, since $C_1$ and $C_2$ are conjugated in $G^+$, and the approach to Chebyshev's bias initiated by Rubinstein and Sarnak \cite{RuSa} and generalized by Ng \cite{Ng}, breaks down. In fact, one can show that $\pi(x, C_1, L/K) = \pi(x, C_2, L/K)$ for any $x \geq 2$ in this case.\\

The random variables at play are symmetric about their mean $m$. This is because their characteristic functions are products of Bessel $J_0$ functions, which are even, multiplied by $t \mapsto \e^{i m t}$ (\cite[Theorem 5.2.1]{Ng}). Therefore, it is easily seen that the corresponding logarithmic density $\delta$ will be $< \frac{1}{2}$, $> \frac{1}{2}$ or $= \frac{1}{2}$ according to whether $m < 0$, $m > 0$ or $m=0$, respectively. They also do not have any atoms, as was shown in \cite[Theorem 2.2]{Dev}.\\

We introduce one more quantity in order to state the main results of \cite{FiJo} on Chebyshev's bias in families of number fields.

\begin{defi} Let $C_1$ and $C_2$ be conjugacy classes of $G$ such that $C_1^+ \neq C_2^+$. Assuming Artin's conjecture, $\GRH$ and $\LIM$, the bias factor of the race between $C_1$ and $C_2$ is $$B(L/K, C_1, C_2) := \frac{\E(X(L/K, C_1, C_2))}{\sqrt{\Var(X(L/K, C_1, C_2))}}.$$ 
\end{defi}

Note that the above variance is non-zero since $C_1^+ \neq C_2^+$.\\

The next theorem is a result giving extremely biased races.

\begin{thm}[\cite{FiJo}, Proposition 5.3]\label{DevFJ} Assume $\GRH, \LI$ and Artin's Conjecture. Then there exists an absolute constant $c_3 > 0$ such that for any conjugacy classes $C_1, C_2$ of $G$ satisfying $C_1^+ \neq C_2^+$ and $B(L/K, C_1, C_2) > 0$, we have $$1 - \delta(L/K, C_1, C_2) < \exp(-c_3B(L/K, C_1, C_2)^2).$$
\end{thm}

\begin{rmq} This theorem is based on a large deviation result due to Montgomery and Odlyzko (see Theorem \ref{DevMO} below). The idea is that if the quantity $B$ is positive and large, then the mean of the associated random variable $X$ is large compared to its variance, and the distribution of the random variable $X$ is very concentrated around its (positive) mean, and therefore it takes positive values with high probability.\\
\end{rmq}

We are actually able to provide lower bounds in the context of Theorem \ref{DevFJ}, but not as uniform as the previous upper bounds (see Theorem \ref{Dev}). This new bound will allow us to exhibit a monotonicity phenomenon in the values of Chebyshev's bias in certain towers of number field extensions (see Theorem \ref{Mono}).\\

The second result can be seen as a converse to the above theorem. We can roughly state it as saying that if the quantity $B$ is small, then the race will be moderately biased, \textit{i.e.} the logarithmic density $\delta$ will be close to $\frac{1}{2}$.

\begin{thm}[\cite{FiJo}, Theorem 5.10]\label{TCL} Assume $\GRH, \LIM$ and Artin's Conjecture. Then for any conjugacy classes $C_1, C_2$ of $G$ such that $C_1^+ \neq C_2^+$, if $|B(L/K, C_1, C_2)|$ is small enough then we have $$\delta(L/K, C_1, C_2) = \frac{1}{2} + \frac{B(L/K, C_1, C_2)}{\sqrt{2\pi}} + O\left(B(L/K, C_1, C_2)^3 + \Var(X(L/K, C_1, C_2))^{-1/3}\right).$$
\end{thm}

\begin{rmqs}\begin{itemize}
\item[i)] This result is a consequence of a central limit behaviour for the random variable $X(L/K, C_1, C_2)$ (explaining the $\sqrt{2\pi}$ factor), provided its variance is large enough, which can be established thanks to bounds on Bessel $J_0$ functions which appear in the characteristic functions of $X(L/K, C_1, C_2)$.
\item[ii)] If $\E(X(L/K, C_1, C_2))$ is bounded and $B(L/K, C_1, C_2)$ approaches zero then the main term besides $\frac{1}{2}$ is $\Var(X(L/K, C_1, C_2))^{-1/3}$. On the other hand, if $\Var(X(L/K; C_1, C_2))^{1/6} = o(|\E(X(L/K; C_1, C_2))|)$ and $\E(X(L/K; C_1, C_2)) = o(\Var(X(L/K; C_1, C_2))^{1/2})$ then we have $\delta(L/K, C_1, C_2) - \frac{1}{2} \sim \frac{B(L/K, C_1, C_2)}{\sqrt{2 \pi}}$ as $B(L/K, C_1, C_2)$ goes to zero.
\end{itemize}
\end{rmqs}

\section{Recollection on ramification and on dihedral and quaternion groups of $2$-power order}

\subsection{Artin conductors}
\label{Conduc}

Let $L/K$ be a Galois extension of number fields with Galois group $G$ and $C_1, C_2$ be distinct conjugacy classes of $G$. In \cite{FiJo}, bounds for the variance of $X(L/K, C_1, C_2)$ are given in terms of the local ramification data of the extension $L/K$. More precisely, a link between this variance and the Artin conductors of the irreducible characters of $G$ is established.

Recall that if $\p$ is a prime ideal of $\O_K$ and $\Pp$ is a prime ideal in $\O_L$ above $\p$, the inertia subgroup $\I(\Pp/\p) := \{g \in G \mid \forall x \in \O_L, g(x)=x \text{ mod } \Pp\} \subset G$ admits a filtration $(G_i(\Pp/\p))_{i \geq 0}$ defined as follows : for integers $i \geq 0$, define $$G_i(\Pp/\p) := \{g \in G \mid \forall x \in \O_L, g(x) = x \text{ mod } \Pp^{i+1}\} \subset G.$$ Obviously, $G_0(\Pp/\p) = \I(\Pp/\p)$ and the subgroups $G_i(\Pp, \p)$ of $G$ only depend on $\Pp$ up to conjugacy, so in the sequel we drop the dependency on $\Pp$ and denote them by $G_i(\p)$. It is also known that the elements of this filtration are eventually trivial. We note that if $L/K$ is tamely ramified, then $G_i(\p)$ is trivial for $i \geq 1$ (\cite[Chapitre IV §2 Corollaire 3]{Ser}). If $\rho : G \to \GL(V)$ is a complex representation, we obtain complex representations of all the ramification subgroups $G_i(\p)$ by restriction. If $\chi$ is the character $\rho$, we define \begin{equation} \label{n} n(\chi, \p) := \sum_{i \geq 0} \frac{|G_i(\p)|}{|G_0(\p)|} \codim(V^{G_i(\p)}),\end{equation} which only depends on $\p$ and not on $\Pp$ since the $G_i(\Pp/\p)$ are conjugates in $G$ and so the dimension of their invariant subspaces are all the same. This sum is finite since the $G_i(\p)$ are eventually trivial, but it is also known that this sum is an integer (\cite[Chapitre VI §2 Corollaire 2]{Ser}), so we can finally define the Artin conductor of the character $\chi$ by $$\mathfrak{f}(L/K, \chi) := \prod_{\p} \p^{n(\chi, \p)}.$$ This ideal of $\O_K$ is well-defined since it is easy to see that for $\p$ unramified in $L$ we have $n(\chi, \p)=0$, so that there are only finitely many prime ideals actually contributing to the product. We recall the important conductor-discriminant formula (\cite[Chapitre VI §3 Corollaire 2]{Ser}).

\begin{thm}\label{Cond} We have \begin{equation} \label{CD} D_{L/K} = \prod_{\chi \in \Irr(G)} \mathfrak{f}(L/K, \chi)^{\chi(1)},\end{equation} where $D_{L/K}$ is the relative discriminant of $L/K$.
\end{thm}

The following quantity appears in the functional equation of the Artin $L$-functions associated to a character $\chi$.

\begin{defi}[\cite{MuMu} p.28] For any irreducible character $\chi$ of $G$, define \begin{equation} \label{A} A(\chi) := |d_K|^{\chi(1)} N_{K/\Q}(\mathfrak{f}(L/K, \chi)),\end{equation} where $d_K$ is the absolute discriminant of $K$.
\end{defi}

\begin{lem}[\cite{FiJo}, Lemma 4.3]\label{B0} Assume Artin's conjecture. Then for any irreducible character $\chi$ of $G$, we have \begin{equation} \label{B_0}B_0(\chi) := \sum_{\gamma_{\chi} \neq 0} \frac{1}{\frac{1}{4} + \gamma_{\chi}^2} \asymp \log A(\chi).\end{equation}
\end{lem}

Recall that the sums $B_0(\chi)$, for non-trivial irreducible characters $\chi$ of $G$, appear in the variances in Theorem \ref{VA}. That explains why we will be interested in bounding the quantities $\log A(\chi)$.

\subsection{Character theory of dihedral and quaternion groups of $2$-power order}

\label{Sect2.2}

We obtain our main contributions in the setting of Galois extensions of number fields with Galois groups dihedral or quaternion of $2$-power order. The goal of this section is to recollect some character-theoretic facts about these groups which will be used throughout Sections 3 and 4. Let us first recall the classical presentations of those groups. 

\begin{defi} Let $n \geq 3$ be an integer. The dihedral group of order $2^n$ is $$D_{2^{n-1}} := \langle r,s \mid r^{2^{n-1}}=s^2=1, srs^{-1}=r^{-1} \rangle.$$ The generalized quaternion group of order $2^n$ is $$\H_{2^n} := \langle x, y \mid x^{2^{n-1}}=1, x^{2^{n-2}} = y^2, yxy^{-1} = x^{-1} \rangle.$$
\end{defi}

The group $\H_8$ is the usual quaternion group of order $8$, with elements usually denoted $i, j, k$ satisfying $$i^2=j^2=k^2=ijk=-1, (-1)^2=1.$$ This notation will be used in Section 3. In $D_{2^{n-1}}$, we will denote by $-1$ the element $r^{2^{n-2}}$, and in $\H_{2^n}$ we will also denote by $-1$ the element $x^{2^{n-2}} = y^2$. Both of those elements are of order $2$ and generate the center of the respective group they belong to.\\

These groups are special instances of metacyclic groups (\textit{i.e.} groups admitting a cyclic normal subgroup with cyclic quotient). As a consequence, computations are easily carried in such groups and the conjugacy classes are easily identified (\cite{Ska}).

\begin{lem}\label{ConjD} Let $n \geq 3$ be an integer. Then $D_{2^{n-1}} = \langle r, s \mid r^{2^{n-1}}=s^2=1, srs^{-1} = r^{-1} \rangle$ has the following $2^{n-2} + 3$ conjugacy classes :

\begin{enumerate}
\item[i)] The trivial conjugacy class $\{1\}$, denoted by $C_1$.
\item[ii)] The conjugacy class $\{-1\}$, denoted by $C_{-1}$.
\item[iii)] Pairs of powers $\{r^k, r^{-k}\}$ for $1 \leq k \leq 2^{n-2}-1$, denoted by $C_{r^k}$.
\item[iv)] The conjugacy class $\{r^ks \mid 0 \leq k \leq 2^{n-2}-1 \text{ even}\}$ of $s$, denoted by $C_s$.
\item[v)] The conjugacy class $\{r^ks \mid 0 \leq k \leq 2^{n-1}-1 \text{ odd}\}$ of $rs$, denoted by $C_{rs}$.
\end{enumerate}
\end{lem} 

\begin{lem}\label{ConjQ} Let $n \geq 3$ be an integer. Then $\H_{2^n} = \langle x, y \mid x^{2^{n-1}}=1, x^{2^{n-2}} = y^2, yxy^{-1} = x^{-1} \rangle$ has the following $2^{n-2} + 3$ conjugacy classes :

\begin{enumerate}
\item[i)] The trivial conjugacy class $\{1\}$, denoted by $C_1$.
\item[ii)] The conjugacy class $\{-1\}$, denoted by $C_{-1}$.
\item[iii)] Pairs of powers $\{x^k, x^{-k}\}$ for $1 \leq k \leq 2^{n-2}-1$, denoted by $C_{x^k}$.
\item[iv)] The conjugacy class $\{x^ky \mid 0 \leq k \leq 2^{n-1}-1 \text{ even}\}$ of $y$, denoted by $C_y$.
\item[v)] The conjugacy class $\{x^ky \mid 0 \leq k \leq 2^{n-1}-1 \text{ odd}\}$ of $xy$, denoted by $C_{xy}$.
\end{enumerate}
\end{lem}

In particular, both $D_{2^{n-1}}$ and $\H_{2^n}$ have $2^{n-2}+3$ isomorphism classes of irreducible complex representations. We refer the reader to \cite{Isa} for classical facts about the representation theory of finite groups. It is a well-known fact that even though those two groups are not isomorphic, they have the same character table.

\begin{lem}\label{Died} Let $n \geq 3$ be an integer and let $D_{2^{n-1}} = \langle r, s \mid r^{2^{n-1}}=s^2=1, srs^{-1} = r^{-1} \rangle$. Let also $\zeta_n = \exp\left(\frac{2i \pi}{2^{n-1}}\right) \in \C$. Then the following homomorphisms are representatives of the $2^{n-2}+3$ isomorphism classes of irreducible complex representations of $D_{2^{n-1}}$ :

\begin{enumerate}
\item[i)] Four abelian representations, coming from the abelianization of $D_{2^{n-1}}$, isomorphic to $\Z/2\Z \times \Z/2\Z$, given as follows : $$\chi_0 : r, s \mapsto 1 \quad \chi_1 : r \mapsto 1, s \mapsto -1 \quad \chi_2 : r \mapsto -1, s \mapsto 1 \quad \chi_3 : r \mapsto -1, s \mapsto -1.$$
\item[ii)] For $1 \leq j \leq 2^{n-2}-1$, a degree $2$ representation given by : $$r \mapsto \begin{pmatrix}\zeta_n^j&0\\0&\zeta_n^{-j}\end{pmatrix}, \quad s \mapsto \begin{pmatrix}1&0\\0&-1\end{pmatrix}$$ with character denoted by $\psi_j$.
\end{enumerate}

The character table of $D_{2^{n-1}}$ is

$$\begin{array}{c|c|c|c|c|c}
 & C_1 & C_{-1} & C_{r^k}, 1 \leq k \leq 2^{n-2}-1 & C_s & C_{rs} \\ 
\hline 
\chi_0 & 1 & 1 & 1 & 1 & 1 \\ 
\hline 
\chi_1 & 1 & 1 & 1 & -1 & -1 \\ 
\hline 
\chi_2 & 1 & 1 & (-1)^k & 1 & -1 \\ 
\hline 
\chi_3 & 1 & 1 & (-1)^k & -1 & 1 \\ 
\hline 
\psi_j, 1 \leq j \leq 2^{n-2}-1 & 2 & (-1)^j2 & \zeta_n^{jk} + \zeta_n^{-jk} = 2\cos\left(\frac{jk\pi}{2^{n-2}}\right) & 0 & 0 \\ 
\end{array}$$
\end{lem}

\begin{lem}\label{Quat} Let $n \geq 3$ be an integer and let $\H_{2^n} = \langle x, y \mid x^{2^{n-1}}=1, x^{2^{n-2}} = y^2, yxy^{-1} = x^{-1} \rangle$. Let also $\zeta_n = \exp\left(\frac{2i \pi}{2^{n-1}}\right) \in \C$. Then the following homomorphisms are representatives of the $2^{n-2}+3$ isomorphism classes of irreducible complex representations of $\H_{2^n}$ :

\begin{enumerate}
\item[i)] Four abelian representations, coming from the abelianization of $\H_{2^n}$, isomorphic to $\Z/2\Z \times \Z/2\Z$, given as follows : $$\chi_0 : x, y \mapsto 1 \quad \chi_1 : x \mapsto 1, y \mapsto -1 \quad \chi_2 : x \mapsto -1, y \mapsto 1 \quad \chi_3 : x \mapsto -1, y \mapsto -1.$$
\item[ii)] For $1 \leq j \leq 2^{n-2}-1$, a degree $2$ representation given by : $$x \mapsto \begin{pmatrix}\zeta_n^j&0\\0&\zeta_n^{-j}\end{pmatrix}, \quad y \mapsto \begin{pmatrix}0&(-1)^j\\1&0\end{pmatrix}$$ with character denoted by $\psi_j$.
\end{enumerate}

The character table of $\H_{2^n}$ is

$$\begin{array}{c|c|c|c|c|c}
 & C_1 & C_{-1} & C_{x^k}, 1 \leq k \leq 2^{n-2}-1 & C_y & C_{xy} \\ 
\hline 
\chi_0 & 1 & 1 & 1 & 1 & 1 \\ 
\hline 
\chi_1 & 1 & 1 & 1 & -1 & -1 \\ 
\hline 
\chi_2 & 1 & 1 & (-1)^k & 1 & -1 \\ 
\hline 
\chi_3 & 1 & 1 & (-1)^k & -1 & 1 \\ 
\hline 
\psi_j, 1 \leq j \leq 2^{n-2}-1 & 2 & (-1)^j 2 & \zeta_n^{jk} + \zeta_n^{-jk} = 2\cos\left(\frac{jk\pi}{2^{n-2}}\right) & 0 & 0 \\ 
\end{array}$$
\end{lem}

We remark that all the irreducible characters of $D_{2^{n-1}}$ and $\H_{2^n}$ are real-valued (which does not necessarily mean they are afforded by real-valued representations, see Lemma \ref{Orth}) and they have bounded degrees (at most $2$) with respect to $n$.

\begin{lem}\label{Fid} With the same notations as in Lemma \ref{Died}, the irreducible character $\psi_j$ of $D_{2^{n-1}}$, for $1 \leq j \leq 2^{n-2}-1$, is faithful if and only if $j$ is odd. The same holds for the irreducible characters of $\H_{2^n}$ as in Lemma \ref{Quat}.
\end{lem}

\begin{demo} Recall that a character $\chi$ of a group $G$ is faithful if the corresponding representation is faithful (that is, injective). This happens if and only if $\chi(g) \neq \chi(1)$ for any $g \in G \setminus \{1\}$.

Obviously, if $j$ is even then $\psi_j(-1)=2=\psi_j(1)$ so $\psi_j$ is not faithful. Conversely, if $j$ is odd, then $\chi(r^k)$ is never equal to $2$, for $1 \leq k \leq 2^{n-2}-1$. Indeed, since $j$ is odd, the $2$-adic valuation of $jk$ is the same as that of $k$, which is $< n-2$, therefore $\frac{jk}{2^{n-1}}$ is not an integer, and $\frac{jk\pi}{2^{n-2}} = \frac{2jk \pi}{2^{n-1}}$ is not an integer multiple of $2\pi$. 
\end{demo}

\begin{lem}\label{Orth} Every irreducible character of $D_{2^{n-1}}$ is orthogonal, and with the same notations as in Lemma \ref{Quat}, the irreducible character $\psi_j$ of $\H_{2^n}$, for $1 \leq j \leq 2^{n-2}-1$, is symplectic if and only if $j$ is odd.
\end{lem}

\begin{demo} Clearly, $\chi_0, \chi_1, \chi_2$ and $\chi_3$ are afforded by real representations so, they are orthogonal (even in the case of $\H_{2^n}$). Consider now a non-abelian character $\psi_j$, $1 \leq j \leq 2^{n-2}-1$, of $D_{2^{n-1}}$. Then by definition $$\varepsilon_2(\psi_j) = \frac{1}{2^n} \sum_{g \in D_{2^{n-1}}} \psi_j(g^2).$$ Since each $r^ks$, with $0 \leq k \leq 2^{n-1}-1$, satisfies $(r^ks)^2=1$ and $1^2=1$, the sum on the right-hand side contains at least $2^{n-1}+1$ terms equal to $\psi_j(1) = 2$. Since this is more than half the number of elements in $D_{2^{n-1}}$ and each $\psi_j(g)$ is smaller than $2$ in absolute value, the reverse triangular inequality implies that $\varepsilon_2(\chi) > 0$. By Theorem \ref{Frob} it must be $1$, and $\psi_j$ is orthogonal.\\

Now let $\psi_j$, $1 \leq j \leq 2^{n-2}-1$, be a non-abelian irreducible character of $\H_{2^n}$. The squares in $\H_{2^n}$ are precisely the even powers of $x$. Indeed, since $y^{-1} = -y$, we have, for any $k \in \{0, \dots, 2^{n-1}-1\}$, $$(x^ky)^2 = -x^kyx^ky^{-1} = -1.$$ From this, and denoting $r(h)$ the number of square roots of $h \in \H_{2^n}$, we see that $r(-1) = 2^{n-1}+2$ while $r(x^{2k})=2$ for $0 \leq k \leq 2^{n-3}-1$. Now if $j$ is odd, then \begin{align*}\varepsilon_2(\psi_j) &= \frac{1}{2^n} \left(4-2(2^{n-1} + 2) + 4\sum_{k=1}^{2^{n-3}-1} \cos\left(\frac{2jk\pi}{2^{n-2}}\right)\right)\\ &= -1 + \frac{1}{2^{n-2}} \sum_{k=1}^{2^{n-3}-1} \cos\left(\frac{2jk\pi}{2^{n-2}}\right).\end{align*} It is easily seen from the triangular inequality that this quantity has to be negative. As before, this shows that $\varepsilon_2(\psi_j)=-1$.

Conversely, if $j$ is even, the reverse triangular inequality implies that $$\varepsilon_2(\psi_j) = \frac{1}{2^n} \left(4+2(2^{n-1} + 2) + 4\sum_{k=1}^{2^{n-3}-1} \cos\left(\frac{2jk\pi}{2^{n-2}}\right)\right) > 0,$$ which implies that $\varepsilon_2(\psi_j)=1$.
\end{demo}

Let us remark that for $n \geq 4$, the group $\H_{2^{n}} = \langle x, y \mid x^{2^{n-1}}=1, x^{2^{n-2}}=y^2, yxy^{-1}=x^{-1} \rangle$ contains $\H_{2^{n-1}}$ as a subgroup generated by $\{x^2, y\}$, and inductively, contains $\H_{2^i} = \langle x^{2^{n-i}}, y \rangle$ for every $3 \leq i \leq n$. Similarly, for $n \geq 3$, $D_{2^{n-2}}$ appears as a subgroup of $D_{2^{n-1}} = \langle r, s \mid r^{2^{n-1}}=s^2=1, srs^{-1} = r^{-1} \rangle$ generated by $\langle r^2, s\rangle$. In particular when $L/K$ is a Galois extension of number field extension with group $\H_{2^n}$ or $D_{2^{n-1}}$, it contains subextensions with Galois groups $\H_{2^i}$ or $D_{2^{i-1}}$ for every $3 \leq i \leq n$. To study Chebyshev's bias in such subextensions, we will need to know how irreducible characters are induced from such $\H_{2^i}$ to $\H_{2^n}$, and from $D_{2^{i-1}}$ to $D_{2^{n-1}}$.

\begin{lem}\label{Ind} With the same notations as in Lemma \ref{Quat}, denote by $\psi_k^{(i)}$, for $3 \leq i \leq n$, $1 \leq k \leq 2^{i-2} -1$, the character of the subgroup $\H_{2^i}$ of $\H_{2^n}$ associated to the representation $$x^{2^{n-i}} \mapsto \begin{pmatrix}\zeta_i^k&0\\0&\zeta_i^{-k}\end{pmatrix} ,\quad y \mapsto \begin{pmatrix}0&-1\\1&0\end{pmatrix},$$ and denote by $\chi_0^{(i)}, \chi_1^{(i)}, \chi_2^{(i)}$ and $\chi_3^{(i)}$ the four abelian characters of $\H_{2^i}$ corresponding to the notations of Lemma \ref{Quat}. For $k \in \{1, \dots, 2^{i-2}-1\}$, one has $$\Ind_{\H_{2^i}}^{\H_{2^n}} \psi_k^{(i)} = \sum_{\underset{l = \pm k \text{ mod } 2^{i-1}}{1 \leq l \leq 2^{n-2}-1}} \psi_l.$$ Also, $$\Ind_{\H_{2^i}}^{\H_{2^n}} \chi_0^{(i)} = \Ind_{\H_{2^i}}^{\H_{2^n}} \chi_1^{(i)} = \chi_0+\chi_1+\chi_2+\chi_3 + \sum_{\underset{j = 0 \text{ mod } 2^{i-1}}{1 \leq j \leq 2^{n-2}-1}} \psi_j$$ and $$\Ind_{\H_{2^i}}^{\H_{2^n}} \chi_2^{(i)} = \Ind_{\H_{2^i}}^{\H_{2^n}} \chi_3^{(i)} = \sum_{\underset{j = 0 \text{ mod } 2^{i-2}, j \neq 0 \text{ mod } 2^{i-1}}{1 \leq j \leq 2^{n-2}-1}} \psi_j.$$ Those equalities also hold when $\H_{2^n}$ and $\H_{2^i}$ are replaced by $D_{2^{n-1}}$ and $D_{2^{i-1}}$ respectively.
\end{lem}

\begin{demo} Let $1 \leq j \leq 2^{n-2}-1$. We use the Frobenius reciprocity formula and the character table in Lemma \ref{Quat} to compute \begin{align*}
\langle \Ind_{\H_{2^i}}^{\H_{2^n}} \psi_k^{(i)}, \psi_{j} \rangle_{\H_{2^n}} = \langle \psi_k^{(i)}, \psi_{j \mid \H_{2^i}} \rangle_{\H_{2^i}} &= \frac{1}{2^i} \sum_{g \in \H_{2^i}} \psi_k^{(i)}(g) \bar{\psi_j(g)}\\ &= \frac{1}{2^i}\left(2\cdot 2 + (-1)^{k+j} \cdot 4 + 2 \sum_{1 \leq l \leq 2^{i-2}-1} \psi_k^{(i)}(x^{2^{n-i}l}) \bar{\psi_j(x^{2^{n-i}l})} \right)\\ &= \frac{1}{2^i}\left(4(1+(-1)^{k+j})\right.\\ &\quad\left.+ 2 \sum_{1 \leq l \leq 2^{i-2}-1} \left(\zeta_i^{kl} + \zeta_i^{-kl}\right)\left(\zeta_n^{jl2^{n-i}} + \zeta_n^{-jl2^{n-i}}\right)\right).
\end{align*}

Since $\zeta_n^{2^m} = \zeta_{n-m}$ for any $0 \leq m \leq n-1$, we find \begin{align*}
\langle \Ind_{\H_{2^i}}^{\H_{2^n}} \psi_k^{(i)}, \psi_{j} \rangle_{\H_{2^n}} &= \frac{1}{2^i}\left(4(1+(-1)^{k+j}) + 2 \sum_{1 \leq l \leq 2^{i-2}-1} \left(\zeta_i^{kl} + \zeta_i^{-kl}\right)\left(\zeta_i^{jl} + \zeta_i^{-jl}\right)\right)\\ &= \frac{1}{2^i}\left(4(1+(-1)^{k+j}) + 2 \sum_{1 \leq l \leq 2^{i-2}-1} \left(\zeta_i^{(k+j)l} + \zeta_i^{-(k+j)l} + \zeta_i^{(k-j)l} + \zeta_i^{(j-k)l}\right)\right)\\ &= \frac{1}{2^i}\left(4(1+(-1)^{k+j}) + 2 \sum_{\underset{l \neq 2^{i-2}}{1 \leq l \leq 2^{i-1}-1}} \left(\zeta_i^{(k+j)l} + \zeta_i^{(k-j)l}\right)\right).
\end{align*}

Now, if $k+j=0$ mod $2^{i-1}$ then $k - j \neq 0$ mod $2^{i-1}$, otherwise we would get $2k = 0$ mod $2^{i-1}$, \textit{i.e.} $k = 0$ mod $2^{i-2}$, which cannot be because $1 \leq k \leq 2^{i-2}-1$. Therefore, the geometric progressions missing two terms sum to \begin{align*}\sum_{\underset{l \neq 2^{i-2}}{1 \leq l \leq 2^{i-1}-1}} \left(\zeta_i^{(k+j)l} + \zeta_i^{(k-j)l}\right) &= 2^{i-1}-2 + \left(\frac{1-\zeta_i^{(k-j)2^{i-1}}}{1 - \zeta_i^{(k-j)}} - \zeta_i^{(k-j)2^{i-2}} - 1\right)\\ &= 2^{i-1} - 4.\end{align*} Similarly, if $k-j=0$ mod $2^{i-1}$, then $k+j \neq 0$ mod $2^{i-1}$, and we find $$\sum_{\underset{l \neq 2^{i-2}}{1 \leq l \leq 2^{i-1}-1}} \left(\zeta_i^{(k+j)l} + \zeta_i^{(k-j)l}\right) = 2^{i-1}-4.$$ Finally, if $k+j \neq 0$ mod $2^{i-1}$ and $k-j \neq 0$ mod $2^{i-1}$, we find $$\sum_{\underset{l \neq 2^{i-2}}{1 \leq l \leq 2^{i-1}-1}} \left(\zeta_i^{(k+j)l} + \zeta_i^{(k-j)l}\right) = -(-1)^{k+j} \cdot 2 - 2.$$

To sum up, we have found $$\langle \Ind_{\H_{2^i}}^{\H_{2^n}} \psi_k^{(i)}, \psi_{j} \rangle_{\H_{2^n}} = \left\{
\begin{array}{l}
1 \text{ if } k+j = 0 \text{ mod } 2^{i-1} \text{ or } k-j=0 \text{ mod } 2^{i-1}\\
0 \text{ otherwise}
\end{array}\right..$$

We have found $2^{n-i}$ irreducible components of $\Ind_{\H_{2^i}}^{\H_{2^n}} \psi_k^{(i)}$ of degree $2$. Since $\Ind_{\H_{2^i}}^{\H_{2^n}} \psi_k^{(i)}(1) = [\H_{2^n} : \H_{2^i}] \psi_k^{(i)}(1) = 2^{n-i+1}$, these are the only ones.\\

The other induced characters are computed in a similar manner, using the Frobenius reciprocity formula. The fact that $\H_{2^n}$ and $D_{2^{n-1}}$ have the same character table with conjugacy classes indexed similarly with respect to their generators imply that the same computations work in the $D_{2^{n-1}}$ case. 
\end{demo}

\begin{cor}\label{Den} Let $3 \leq i \leq n-1$. If $\chi \in \Irr(\H_{2^i})$ then for any $\chi' \in \Irr(\H_{2^i}) \setminus \{\chi\}$, inducing from $\H_{2^i}$ to $\H_{2^n}$ we have $(\chi, \chi')=1$ unless $\deg(\chi)=1$, where $(\chi, \chi')=1$ means that $\langle \Ind_{\H_{2^i}}^{\H_{2^n}} \chi, \Ind_{\H_{2^i}}^{\H_{2^n}} \chi' \rangle = 0$ (see the paragraph above Theorem \ref{Dev}).  This also holds when $\H_{2^n}$ and $\H_{2^i}$ are replaced by $D_{2^{n-1}}$ and $D_{2^{i-1}}$ respectively.
\end{cor}

The above conditions will allow us to obtain a lower bound for $1 - \delta$ for relative number field extensions (i.e. with base field different from $\Q$) with Galois groups $\H_{2^i}$ and $\mathcal{D}_{2^{i-1}}$ using Theorem \ref{Dev}.\\

The following lemma will be used to compute moments in Section 4.

\begin{lem}\label{Root} For $i \geq 3$ and $1 \leq k \leq 2^{i-2}-1$, $$\sum_{\underset{j \text{ odd}}{j=1}}^{2^{i-2}-1} \left(\zeta_{i}^{jk} + \zeta_{i}^{-jk}\right) = 0.$$
\end{lem}

\begin{demo} This is clear for $i=3$, so assume $i \geq 4$. This sum can be rewritten $$\sum_{j=1}^{2^{i-2}-1} \left(\zeta_{i}^{jk} + \zeta_{i}^{-jk}\right) - \sum_{h=1}^{2^{i-3}-1} \left(\zeta_{i}^{2hk} + \zeta_{i}^{-2hk}\right).$$ As in the previous proof, we see that the first sum is simply $$\sum_{j=0}^{2^{i-1}-1} \zeta_{i}^{jk} - \zeta_{i}^0 - \zeta_{i}^{k2^{i-2}} = -1 - (-1)^k$$ since $0< k < 2^{i-1}$ and so $\zeta_{i}^k \neq 1$. The second sum $$\sum_{h=1}^{2^{i-3}-1} \left(\zeta_{i}^{2hk} + \zeta_{i}^{-2hk}\right)$$ is being dealt with similarly, because $\zeta_{i}^2 = \zeta_{i-1}$ and $0 < k < 2^{i-2}$, so it also equals $-1-(-1)^k$. Thus the two sums cancel each other.
\end{demo}

\section{Extensions of $\Q$ of group $\H_8$ : horizontal Chebotarev biases}

\label{Sect3}

In this section, we will depart temporarily from the notations we used in Section 2, and we will use the more familiar $i, j, k$ notations for elements of $\H_8$. Our goal is to prove Theorem \hyperlink{ThA}{A}. We will be constructing two distinct families of number fields with Galois group over $\Q$ isomorphic to $\H_8$ and with opposite biases along each family, due to the existence or not of a central zero for the corresponding Dedekind zeta functions.\\

Recall that $\mathbb{H}_8 = \{\pm 1, \pm i, \pm j, \pm k\}$ is the quaternion group of order $8$. Its elements all commute with $-1$ and satisfy the following relations $$i^2=j^2=k^2=-1,\, ij=k,\, ji=-k,\, jk=i,\, kj=-i,\, ki=j \text{ and } ik=-j.$$

Here is the lattice of subgroups of $\H_8$ (see \cite{Deb})

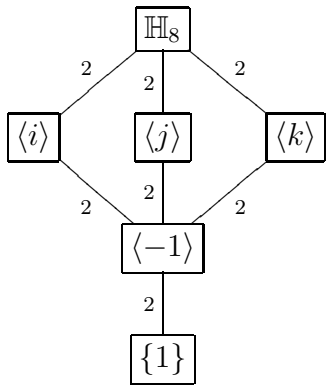
$$\xymatrix{
&*+[F]{\mathbb H_8} \\
*+[F]{\langle i \rangle} \ar@{-}[ru]^2 & *+[F]{\langle j \rangle} \ar@{-}[u]^2 & *+[F]{\langle k \rangle} \ar@{-}[lu]_2 \\
&*+[F]{\langle -1 \rangle} \ar@{-}[ru]_2 \ar@{-}[u]^2 \ar@{-}[lu]^2 \\
&*+[F]{\{1\}} \ar@{-}[u]^2 
}$$

\noindent The character table of $\H_8$ is

$$\begin{tabular}{c|c|c|c|c|c}
 & $C_1$ & $C_{-1}$ & $C_i$ & $C_j$ & $C_k$ \\ 
\hline 
$\chi_0$ & 1 & 1 & 1 & 1 & 1 \\ 
\hline 
$\chi_i$ & 1 & 1 & 1 & -1 & -1 \\ 
\hline 
$\chi_j$ & 1 & 1 & -1 & 1 & -1 \\ 
\hline 
$\chi_k$ & 1 & 1 & -1 & -1 & 1 \\ 
\hline 
$\psi$ & 2 & -2 & 0 & 0 & 0 \\ 
\end{tabular} $$

\noindent where $C_x$ denotes the conjugacy class of $x \in \mathbb H_8$.\\

Let $K/\Q$ be a Galois extension, with Galois group $G$ isomorphic to $\H_8$. Assuming $\LI$, we can compute the values of the different means and variances, according to the formulas given in Theorem \ref{VA} (since $$\E(X(K/\Q, C_1, C_2)) = -\E(X(K/\Q, C_2, C_1)) \text{ and } \Var(X(K/\Q, C_1, C_2))=\Var(X(K/\Q, C_2, C_1)),$$ there are only $\binom{5}{2}=10$ races to consider).

\begin{prop} Assume $\LI$ for the extension $K/\Q$. Let $\z = \ord_{s= \frac{1}{2}} L(s, \psi, K/\Q)$, where $\psi$ is the non-abelian character of $\H_8$ in the above character table. Then one has

$$\begin{array}{c|c|c|c}
a&b&\E(X(K/\Q, C_a, C_b))&\Var(X(K/\Q, C_a, C_b))\\
\hline
1&-1& 4(1-2\z)&16B_0(\psi)\\
\hline
1&i,j,k&-2(1+2\z)&4\sum_{\chi \neq \chi_b} B_0(\chi)\\
\hline
-1&i,j,k&2(2\z-3)&4\sum_{\chi \neq \chi_b} B_0(\chi)\\
\hline
i,j,k&i,j,k&0&4B_0(\chi_a) + 4B_0(\chi_b)
\end{array}$$

where $B_0$ is defined in \ref{B_0}.
\end{prop}

\begin{rmqs} Let us make a few comments on these values : \begin{itemize}
\item[i)] The presence or not of a zero at $\frac{1}{2}$ for $s \mapsto L(s, \psi, K/\Q)$ changes the sign of $\E(X(K/\Q, C_1, C_{-1}))$.
\item[ii)] The race between $C_1$ and $C_b$, for $b \in \{i, j, k\}$, is always biased towards $C_b$. This was expected since $1$ is a square in $G$, and $b$ is not. 
\item[iii)] The presence or not of a zero  with multiplicity at least $2$ at $\frac{1}{2}$ for $s \mapsto L(s, \psi, K/\Q)$ changes the sign of $\E(X(K/\Q, C_{-1}, C_b))$, for $b \in \{i, j, k\}$.
\item[iv)] There is no bias in the race between $C_a$ and $C_b$, for $a \in \{i, j, k\}$ and $b \in \{i, j, k\} \setminus \{a\}$.\\
\end{itemize}
\end{rmqs}

Since we will be using the hypothesis $\LIP$, under which the order of vanishing $\z$ can only be $0$ or $1$, and we want to exhibit a change of leading conjugacy class in Chebotarev races, we will now focus on the conjugacy classes $C_1$ and $C_{-1}$. When $K/\Q$ is tamely ramified, we can deduce more precise bounds for the variances.

\begin{prop}\label{Tab8} Assume $\LIM$ and that $K/\Q$ is tamely ramified. Then $$\Var(X(K/\Q, C_1, C_{-1})) \asymp \log |d_K|.$$
\end{prop}

\begin{demo} Since $\K/\Q$ tame, the filtration of the inertia subgroup (see Section 2.1) at any prime number ramified in $K$ only has length $1$, \textit{i.e.} for any prime $p$ and $i \geq 1$, $|G_i(p)|=1$. In particular, for any non-trivial character $\chi$ of $G$ and any prime number $p$, one has $$n(\chi, p) = \codim V^{\I(p)},$$ where $V$ is the space of the representation affording $\chi$.

Using the definition, we find for any prime $p$ ramified in $K$ (so that $\I(p) = G_0(p) \neq \{1\}$), $$n(\chi_i, p) = \left\{ \begin{array}{l}
1 \text{ if } \I(p) = \langle j \rangle \text{ or } \I(p) = \langle k \rangle\\
0 \text{ otherwise},
\end{array}\right.$$ $$n(\chi_j, p) = \left\{ \begin{array}{l}
1 \text{ if } \I(p) = \langle i \rangle \text{ or } \I(p) = \langle k \rangle\\
0 \text{ otherwise},
\end{array}\right.$$ $$n(\chi_k, p) = \left\{ \begin{array}{l}
1 \text{ if } \I(p) = \langle i \rangle \text{ or } \I(p) = \langle j \rangle\\
0 \text{ otherwise}
\end{array}\right.$$ and $$n(\psi, p) = 2.$$

This yields $$\mathfrak{f}(K/\Q, \chi_i) = \prod_{\underset{\I(p) = \langle j \rangle}{p \mid d_K}} p \times \prod_{\underset{\I(p) = \langle k \rangle}{p \mid d_K}} p,$$ $$\mathfrak{f}(K/\Q, \chi_j) = \prod_{\underset{\I(p) = \langle i \rangle}{p \mid d_K}} p \times \prod_{\underset{\I(p) = \langle k \rangle}{p \mid d_K}} p,$$ $$\mathfrak{f}(K/\Q, \chi_k) = \prod_{\underset{\I(p) = \langle i \rangle}{p \mid d_K}} p \times \prod_{\underset{\I(p) = \langle j \rangle}{p \mid d_K}} p$$ and $$\mathfrak{f}(K/\Q, \psi) = \prod_{p \mid d_K} p^2.$$

Since the base field is $\Q$, we have $A(\chi)=\mathfrak{f}(K/\Q, \chi)$ for any character $\chi$ of $G$ by (\ref{A}), and the conductor-discriminant formula (\ref{CD}) thus gives $$|d_K| = \prod_{\underset{e_p = 8}{p \mid d_K}} p^4 \prod_{\underset{e_p = 4}{p \mid d_K}} p^6 \prod_{\underset{e_p = 2}{p \mid d_K}} p^4,$$ where $e_p$ is the ramification index of the ramified prime $p$. In particular, we have $$A(\psi)^2 \leq d_K \leq A(\psi)^3,$$ so that \begin{equation}\label{log}\log A(\psi) \asymp \log d_K.\end{equation} Finally, Theorem \ref{VA} and Lemma \ref{B0} yield $$\Var(X(K/\Q, C_1, C_{-1})) \asymp B_0(\psi) \asymp \log A(\psi) \asymp \log |d_K|.$$
\end{demo}

The fact that a bound such as Proposition \ref{Tab8} holds follows from the fact that the character $\psi$ is faithful, because this implies that each ramified prime appears with positive exponent in $A(\chi) = \mathcal f(K/\Q, \psi)$. This observation will be used again in our study of biases in towers of extensions (see Proposition \ref{VarQ}).\\

We will use the following theorem of Fröhlich to construct the quaternion extensions of $\Q$ we will be interested in :

\begin{thm}[\cite{Fro1}]\label{Fro1} Let $d$ be a square-free integer, $d > 1$ and $d = 1$ mod $4$. Let $\mathfrak{R}$ be a finite set of primes unramified in $\Q(\sqrt d)$ and not containing $2$. Then there exist infinitely many tamely ramified extensions $K/\Q$ and $L/\Q$ such that :

\begin{itemize}
\item[i)] $\Gal(K/\Q) \simeq \Gal(L/\Q) \simeq \H_8$.
\item[ii)] $\Q(\sqrt{d}) \subset K \cap L$.
\item[iii)] Every prime in $\mathfrak{R}$ is ramified in $K$ and $L$.
\item[iv)] $W(\psi, K/\Q)=-1$ and $W(\psi, L/\Q)=1$.
\end{itemize}
\end{thm}

\begin{rmq} Theorem \ref{Fro1} is in fact a weak version of Fröhlich's theorem in \cite{Fro1}. Actually, one can specify any finite number of unramified primes in $K$ and $L$ as well (the ramification must be compatible with the fact that $\Q(\sqrt{d}) \subset K \cap L$), and one can also ask for $K$ and $L$ to be totally real or totally imaginary.
\end{rmq}

\begin{cor}\label{Corps} For any square-free integer $d > 1$ with $d=1$ mod $4$ and any prime $p$, there exist infinitely many tamely ramified extensions $K/\Q$ and $L/\Q$ such that :

\begin{itemize}
\item[i)] $\Gal(K/\Q) \simeq \Gal(L/\Q) \simeq \H_8$.
\item[ii)] $\Q(\sqrt{d}) \subset K \cap L$.
\item[iii)] $p$ is ramified in $K$ and $L$.
\item[iv)] $L\left(\frac{1}{2}, \psi, K/\Q\right)=0$ and if $\LIP$ holds then $L\left(\frac{1}{2}, \psi, L/\Q\right) \neq 0$.
\end{itemize}
\end{cor}

\begin{demo} The existence of the number fields $K$ follows from Theorem \ref{Fro1} without the assumption of $\LIP$ : we consider tamely ramified extensions $K/\Q$ with Galois group $\H_8$, ramified at $p$, containing $\Q(\sqrt d)$ and such that $W(\psi)=-1$. Proposition \ref{Zerr} shows that $L(1/2, \psi, K/\Q)=0$.

The existence of the number fields $L$ follows similarly because, assuming $\LIP$, $W(\psi)=1$ implies that $s \mapsto L(s, \psi, L/\Q)$ does not vanish at $\frac{1}{2}$.
\end{demo}

We are now ready to prove Theorem \hyperlink{ThA}{A}, which we recall for convenience.

\begin{thm}\label{Opp8} Assume $\GRH$ and $\LIP$. For any function $f$ such that $f(n) \underset{n \to +\infty}{\longrightarrow} +\infty$, there exist two infinite families $(K_d)_d$ and $(L_d)_d$ of Galois extensions of $\Q$, indexed by square-free integers satisfying $d > 1$ and $d=1$ mod $4$, such that for any $d$ in the index set :

\begin{itemize}
\item[i)] $\Q(\sqrt d) \subset K_d \cap L_d$.
\item[ii)] $\Gal(K_d/\Q) \simeq \Gal(L_d/\Q) \simeq \H_8$.
\item[iii)] $0 < \frac{1}{2} - \delta(K_d/\Q, C_1, C_{-1}) \ll \frac{1}{f(d)}$.
\item[iv)] $0 < \delta(L_d/\Q, C_1, C_{-1}) - \frac{1}{2} \ll \frac{1}{f(d)}$.
\end{itemize}
\end{thm}

\begin{demo} For any square-free $d > 1$ with $d=1$ mod $4$, we choose a field $K_d$ given in the first part of Corollary \ref{Corps}, ramified at the smallest prime larger than $e^{f(d)^3}$. By construction of $K_d$ we have $$B(X(K_d/\Q, C_1, C_{-1})) = \frac{\E(X(K_d/\Q, C_1, C_{-1}))}{\sqrt{\Var(X(K_d/\Q, C_1, C_{-1}))}} = \frac{1-2\z}{B_0(\psi)^{1/2}} = - \frac{1}{B_0(\psi)^{1/2}}.$$ In particular, $\E(X(K_d/\Q, C_1, C_{-1}))$ is negative, because of the existence of a real zero of $s \mapsto L(s, \psi, K_d/\Q)$, and this implies that $\delta(K_d/\Q, C_1, C_{-1}) < \frac{1}{2}$. Lemma \ref{B0} yields $$B_0(\psi) \asymp \log A(\psi),$$ and by (\ref{log}) we have $$\log A(\psi) \asymp \log |d_{K_d}|,$$ since $\psi$ is faithful and $K_d/\Q$ is tamely ramified. Moreover, since $K_d/\Q$ is ramified at the smallest prime larger than $e^{f(d)^3}$, we have $$\log |d_{K_d}| \gg \log\left(e^{f(d)^3}\right) = f(d)^3.$$ Therefore, $$|B(X(K_d/\Q, C_1, C_{-1}))| \ll \frac{1}{f(d)^{3/2}}.$$ Similarly, we have $$\Var(X(K_d/\Q, C_1, C_{-1}))^{-1} \ll \frac{1}{f(d)^3}.$$ Since Artin's conjecture holds for extensions with Galois group $\H_8$ (recall that this group is supersolvable, and that Artin's conjecture is known in that case), Theorem \ref{TCL} implies that $$0 < \frac{1}{2} - \delta(K_d/\Q, C_1, C_{-1}) \ll \frac{1}{f(d)}.$$

The construction of $L_d$ follows along the same lines, by choosing fields as in the second part of Corollary \ref{Corps}. In that case, the mean is positive because there is no real zero of $s \mapsto L(s, \psi, L_d/\Q)$, so $\delta(K_d/\Q, C_1, C_{-1}) > \frac{1}{2}$, and the estimates are the same.
\end{demo}

\begin{rmq} If we had a statement analoguous to Theorem \ref{Fro1} in which we could specify if $s \mapsto L(s, \psi, K/\Q)$ vanishes, or not, at $\frac{1}{2}$ with multiplicty at least two, then we would be able to obtain a result similar to Theorem \hyperlink{ThA}{A} for the race between $C_{-1}$ and $C_b$, for $b \in \{i, j, k\}$, \textit{i.e.} moderately biased races with biases of opposite signs, bounded in absolute value by $\frac{1}{f(d)},$ for any function $f$ going to infinity. Such a result would of course contradict $\LIP$.
\end{rmq}

\section{Chebyshev's bias in towers}

We now consider the vertical aspect of our problem. Our goal is to prove Theorems \hyperlink{ThB}{B}, \hyperlink{ThC}{C} and \hyperlink{ThD}{D}. Instead of working with a family of extensions of $\Q$ with fixed Galois group, we want to study different Galois extensions of $\Q$ with Galois groups dihedral of order a power of two or generalized quaternion of increasing sizes. We also want to compare Chebyshev's bias in subextensions. For this to make sense, we need to make a few observations first.

\begin{lem} For $n \geq 4$, if we let $\H_{2^{n}} = \langle x, y \mid x^{2^{n-1}}=1, x^{2^{n-2}}=y^2, yxy^{-1} = x^{-1}\rangle$, then $\langle x^2, x^ky \rangle$ is a normal subgroup of $\H_{2^{n}}$ isomorphic to $\H_{2^{n-1}}$, for any $0 \leq k \leq 2^{n-2}-1$.
\end{lem} 

\begin{demo} It is easily seen that $(x^2)^{2^{n-2}} = 1, (x^ky)^2 = x^kyx^ky^{-1}y^2 = y^2 = (x^2)^{2^{n-3}}$ and $x^ky x^2 (x^ky)^{-1} = x^kyxy^{-1} y x y^{-1} x^{-k} = x^{-2}$. The map $x^ay^b \mapsto a + kb \text{ mod } 2$ is a group homomorphism whose kernel is $\langle x^2, x^ky \rangle$, which proves that it is normal in $\H_{2^n}$ and of order $2^{n-1}$. Since it satisfies the same relations as $\H_{2^{n-1}}$ and has the same order, it is isomorphic to it.
\end{demo}

Therefore, if $K/\Q$ is a number field extension with $\Gal(K/\Q) \simeq \H_{2^n}$ then there exist number fields $K \supset K_3 \supset K_4 \supset \dots \supset K_{n-1} \supset K_n = \Q$ with $\Gal(K/K_i) \simeq \H_{2^i}$ for all $3 \leq i \leq n$, namely the subfields fixed by $\langle x^{2^{n-i}}, y \rangle$. In the following, we will always assume that the quaternion subgroups have been chosen so that $\Gal(K/K_i)$ is generated by $x^{2^{n-i}}$ and $y$, for $3 \leq i \leq n-1$.
Similarly, if $n \geq 3$ then $D_{2^{n-1}}$ is a normal subgroup of $D_{2^n} = \langle r,s \mid r^{2^{n-1}}=s^2, srs^{-1}=r^{-1} \rangle$, generated for example by $r^2$ and $s$. Those observations will allow us to compare Chebyshev's bias in subextensions (see Theorem \ref{Mono}).

\subsection{Moments in generalized quaternion Galois groups}

\label{Sect4.1}

In this section, we compute bounds on the moments of the random variables attached to the different Chebotarev races in number field extensions with Galois group generalized quaternion. They will be used in the proof of Theorem \hyperlink{ThC}{C}.\\

We will use the following fact about root numbers of symplectic characters of generalized quaternion groups.

\begin{thm}[\cite{Fro2}, Theorem 3]\label{Fro2} Let $N/L$ be a tamely ramified extension of number fields, with Galois group generalized quaternion. Then the root numbers of the symplectic irreducible characters of $\Gal(N/L)$ are all equal.
\end{thm}

Now, let $n \geq 3$, and assume that $K/\Q$ is a Galois extension with Galois group $$G := \langle x, y \mid x^{2^{n-1}}=1, x^{2^{n-2}} = y^2, yxy^{-1} = x^{-1} \rangle \simeq \H_{2^n},$$ and let $\Q = K_n \subset \dots \subset K_3 \subset K$ be intermediate number fields as before, that is $$G_i := \langle x^{2^{n-i}}, y \rangle \simeq \H_{2^i}$$ and $$K_i := K^{G_i},$$ for $3 \leq i \leq n$. We will denote $x^{2^{n-i}}$ by $x_i$ so that $G_i = \langle x_i, y\rangle.$ Conjugacy classes in $G_i$ will be denoted, according to our notations in Lemma \ref{ConjQ}, by $C_1^{(i)}, C_{-1}^{(i)}, C_{x_i^k}^{(i)}, C_y^{(i)}$ and $C_{x_i y}^{(i)}$. Note that each $G_i$ is supersolvable, so Artin's conjecture is known to hold for Artin $L$-functions attached to irreducible characters of $G_i$.\\

Recall that if $C$ is a conjugacy class of $G_i$, then $C^+$ is the conjugacy class $\bigcup_{g \in G} gCg^{-1}$ of $G_i^+=G$. The following lemma is immediate.

\begin{lem}\label{C+} Let $i \in \{3, \dots, n-1\}$. Using the notations of section 2, we have \begin{itemize}
\item[i)] $C_1^{(i) +} = C_1$.
\item[ii)] $C_{-1}^{(i)+} = C_{-1}$.
\item[iii)] For any $1 \leq k \leq 2^{i-2}-1$, $C_{x_i^k}^{(i)+} = C_{x_i^k}$.
\item[iv)] $C_y^{(i)+} = C_y$.
\item[v)] $C_{x_iy}^{(i)+} = C_{y}$.
\end{itemize}

In particular, if $C_1$ and $C_2$ are distinct conjugacy classes of $G_i$, then $C_1^+ \neq C_2^+$, unless, $\{C_1, C_2\} = \{C_y^{(i)}, C_{x_iy}^{(i)}\}$. 
\end{lem}

When $K/\Q$ is tamely ramified, we denote by $W_{K/K_i}$ the root number of any of the symplectic characters of $\Gal(K/K_i)$. We first relate the root numbers of symplectic characters of $G$ to those of each $G_i$, and the orders of vanishing at $1/2$ of the corresponding Artin $L$-functions.

\begin{prop}\label{Zero} Assume $\LIP$ and assume $K/\Q$ is tamely ramified. Then for any $3 \leq i \leq n$, we have $W_{K/K_i} = W_K$, where $W_K = W_{K/\Q}$ is the root number of any symplectic character of $G$. Moreover, if $W_K=-1$ then only the symplectic characters of $G_i$ have their Artin $L$-function vanish at $1/2$, and the order of vanishing is $2^{n-i}$.
\end{prop}

\begin{demo} Clearly, $K/K_i$ is tamely ramified so our definition of $W_{K/K_i}$ makes sense using Theorem \ref{Fro2}. For any $3 \leq j \leq n$, consider the classical decomposition \begin{align*}\zeta_{K}(s) &= L(s, \chi_0, K/K)\\ &= L(s, \Ind_{\{1\}}^{G_j} \chi_0, K/K_j)\\ &= L(s, \reg_{G_j}, K/K_j)\\ &= L(s, \sum_{\chi \in \Irr(G_j)} \chi(1) \chi, K/K_j)\\ &= \prod_{\chi \in \Irr(G_j)} L(s, \chi, K/K_j)^{\chi(1)},
\end{align*} where $\reg_{G_j}$ is the character of the regular representation of $G_j$.

If $W_K=1$, consider the previous factorization with $j=n$. By $\LIP$, none of the factors vanish at $1/2$, and so $\zeta_K$ does not vanish at $1/2$. Using now the decomposition with $j=i$, we see that $L(1/2, \chi, K/K_i) \neq 0$ for each irreducible character $\chi$ of $G_i$. In particular, $W_{K/K_i} \neq -1$, i.e. $W_{K/K_i}=1$.

Conversely, if $W_K=-1$, then $\zeta_K$ vanishes at $1/2$ to order $2^{n-2}$. Indeed, under $\LIP$, the only factors that vanish at $1/2$ are the $L(s, \chi, K/\Q)^{\chi(1)}$ for $\chi=\chi_0$ or $\chi$ symplectic. There are $2^{n-3}$ symplectic characters, all vanishing at $1/2$ to order one (again, because of $\LIP$, and the fact all such characters have root number $-1$) and satisfying $\chi(1)=2$. Moreover, $L(s, \chi_0, K/\Q) = \zeta(s)$ is the classical Riemann $\zeta$ function, which is known not to vanish at $1/2$. Using Lemma \ref{Ind}, we see that for $1 \leq k \leq 2^{i-2}-1$ odd we have \begin{align*}L(s, \psi_k^{(i)}, K/K_i) &= L(s, \Ind_{G_i}^G \psi_k^{(i)}, K/\Q)\\ &= L(s, \sum_{\underset{l = \pm k \text{ mod } 2^{i-1}}{0 \leq l \leq 2^{n-2}-1}} \psi_l, K/\Q)\\ &= \prod_{\underset{l = \pm k \text{ mod } 2^{i-1}}{0 \leq l \leq 2^{n-2}-1}} L(s, \psi_l, K/\Q).
\end{align*} Since each $L(s, \psi_l, K/\Q)$, with $1 \leq l \leq 2^{n-2}-1$ odd, vanishes at $1/2$ to order one, we see that  each $L(s, \psi_k^{(i)}, K/K_i)$ vanishes at $1/2$ as well, to order $2^{n-i}$, for $1 \leq k \leq 2^{i-2} -1$ odd. In particular, $W_{K/K_i}=-1$. Moreover, the $2^{i-3}$ symplectic characters of $G_i$ contribute to the vanishing of $\zeta_K$ at $1/2$ to order $2 \times 2^{i-3} \times 2^{n-i} = 2^{n-2}$, so no other irreducible character of $G_i$ have its Artin $L$-function vanishing at $1/2$.
\end{demo}

In Propositions \ref{VarQ}, \ref{Esp1} and \ref{EspQ}, we give bounds on the variances and compute the means of the random variables attached to the Chebotarev races in the extensions $K/K_i$.

\begin{prop}\label{VarQ} Assume $\GRH$ and $\LIM$. Then for any $3 \leq i \leq n$ and any distinct conjugacy classes $C_1$, $C_2$ of $G_i$ such that $\{C_1, C_2\} \neq \{C^{(i)}_{y}, C^{(i)}_{x_i y}\}$, $$\Var(X(K/K_i, C_1, C_2)) \ll \log |d_K|.$$ 
Moreover if $K/\Q$ is tamely ramified and $\{C_1, C_2\} \not \subset \{C^{(i)}_{x_i^{2^{i-3}}}, C^{(i)}_{y}, C^{(i)}_{x_i y}\}$ then for any $3 \leq i \leq n,$ $$\Var(X(K/K_i, C_1, C_2)) \gg \frac{\log |d_K|}{16^{n}},$$ and if $\{C_1, C_2\} \cap \{C_{x_i^k}^{(i)} \mid 1 \leq k \leq 2^{i-2}-1\} = \emptyset$ this can be improved to $$\Var(X(K/K_i, C_1, C_2)) \gg \log |d_K|.$$
\end{prop}

\begin{demo} By Lemma \ref{C+}, the condition $\{C_1, C_2\} \neq \{C^{(i)}_{y}, C^{(i)}_{x_i y}\}$ ensures that $C_1^+ \neq C_2^+$. Recall that $$\Var(X(K/K_i, C_1, C_2)) = \sum_{\chi \in \Irr(G^+)} |\chi(C_1^+)-\chi(C_2^+)|^2 B_0(\chi)$$ by Theorem \ref{VA}, with $$B_0(\chi) \asymp \log A(\chi)$$ by Lemma \ref{B0}. In particular, we see using Lemma \ref{C+} that our estimates on $\Var(X(K/K_i, C_1, C_2))$ do not depend on $K_i$, so it is enough to prove them in the case $i=n$, that is $K_i=\Q$.

Since the values taken by the irreducible characters of $G$ are bounded in absolute value uniformly in $n$, we get $$\Var(X(K/\Q, C_1, C_2)) \ll \sum_{\chi \in \Irr(G)} \log A(\chi) \leq \sum_{\chi \in \Irr(G)} \chi(1) \log A(\chi).$$ Recall that, according to (\ref{A}) for any irreducible character $\chi$ of $G$, $$A(\chi) = \mathfrak{f}(K/\Q, \chi)$$ is the Artin conductor of $\chi$, so the conductor-discriminant formula (\ref{CD}) yields \begin{align*}\sum_{\chi \in \Irr(G)} \chi(1) \log A(\chi) &= \sum_{\chi} \chi(1) \log \mathfrak{f}(K/\Q, \chi)\\ &= \log |d_K|.\end{align*} This proves the stated upper bound.

To prove the first lower bound, assume $K/\Q$ is tamely ramified. If $\chi$ is an irreducible symplectic character of $G$ then $\chi$ is faithful by Lemma \ref{Fid}. This implies that there are no invariant vectors for the representation  of character $\chi$. Moreover, since $K/\Q$ is tamely ramified, the ramification groups of index $\geq 2$ are trivial for any prime $p$ ramified in $K$, and for any such prime we find $n(\chi, p) = 2$ (where $n(\chi, p)$ is defined in (\ref{n})). In particular, $$A(\chi) = \mathfrak{f}(K/\Q, \chi) = \prod_{p \mid d_K} p^2.$$ Now, since $K/\Q$ is tamely ramified, each prime $\p$ of $K$ above a ramified prime number $p$ appears in the factorization of $\partial_{K/\Q}$, the different of $K/\Q$, with exponent $e_{\p}-1$ where $e_{\p}$ denotes the corresponding ramification index (\cite[Theorem 2.6]{Neu}). Since the discriminant of $K/\Q$ is the $K/\Q$-norm of the different, this yields $$|d_{K}| = \prod_{\p \mid \partial_{K/\Q}} N_{K/\Q}(\p)^{e_{\p}-1} = \prod_{p \mid d_K} p^{(e_{p}-1)f_{p}g_{p}},$$ where $f_{p}$ (respectively $g_{p}$) denotes the residual degree of (respectively the number of primes above) the prime $p$. In particular, this exponent is less than $e_{p}f_{p}g_{p} = [K:\Q]=2^n$. Therefore we have \begin{align*}\prod_{p \mid d_K} p^{2^n} &\geq \prod_{p \mid d_K} p^{(e_p-1)f_pg_p}\\&= |d_K| .\end{align*} Recall there are exactly $2^{n-3}$ symplectic characters of $G$, so we find \begin{align*}\prod_{\chi \text{ symplectic}} A(\chi) &= \prod_{p \mid d_K} p^{2^{n-2}}\\ &\geq |d_K|^{1/4},\end{align*}
and we find $$\sum_{\chi \text{ symplectic}} \log A(\chi) \gg \log |d_K|.$$

We now need to lower bound the quantity $|\chi(C_1)-\chi(C_2)|$ as $\chi$ ranges over the set of symplectic characters of $G$. A straightforward inspection of all cases shows that $\chi(C_1) = \chi(C_2)$ can only happen if $\{C_1, C_2\} \subset \{C_{x^{2^{n-3}}}, C_y, C_{xy}\}$. In all the other cases, we see that the minimum value of $|\chi(C_1)-\chi(C_2)|$, when $\chi$ varies in the set of symplectic irreducible characters of $G$, is $\gg \frac{1}{4^{n}}$. Indeed, if $\{C_1, C_2\} = \{C_{x^k}, C_{x^l}\}$ with $1 \leq k < l \leq 2^{n-2}-1$, then by the mean value theorem, we have for any $1 \leq j \leq 2^{n-2}-1$, $$|\psi_j(C_1) - \psi_j(C_2)| = 2\left|\cos\left(\frac{jk \pi}{2^{n-2}}\right) - \cos\left(\frac{jl \pi}{2^{n-2}}\right)\right| = \frac{j|k-l|\pi}{2^{n-3}} |\sin(\xi_n)|$$ for some $\frac{jk \pi}{2^{n-2}} < \xi_n < \frac{jl\pi}{2^{n-2}}$. By $\pi$-periodicity of $|\sin|$ and the fact that $x \mapsto \left|\sin\left(\frac{\pi}{2} - x\right)\right|$ is even, it is easy to see that $|\sin(\xi_n)| = |\sin(\xi_n')|$ for some $\frac{\pi}{2^{n-2}} < |\xi_n'| < \frac{\pi}{2}$ (\textit{i.e.} we move to the first quadrant). By the classical inequality $|\sin(x)| \geq \frac{2}{\pi} |x|$ for $|x| \leq \frac{\pi}{2}$ and the fact that $|k-l| \geq 1$ and $j \geq 1$ we finally get $|\psi_j(C_1) - \psi_j(C_2)| \gg \frac{1}{4^n}$. The cases involving $C_1, C_{-1}, C_y$ and $C_{xy}$ are similar. This yields $$\Var(X(K/\Q, C_1, C_2)) \gg \frac{1}{16^{n}} \sum_{\chi \text{ symplectic}} \log A(\chi) \gg \frac{\log |d_K|}{16^n}.$$ As for the last lower bound, if $\{C_1, C_2\} \cap \{C_{x^k} \mid 1 \leq k \leq 2^{n-2}-1\} = \emptyset$ then we actually have $|\chi(C_1)-\chi(C_2)| \geq 2$ for $\chi$ irreducible symplectic, so $$\Var(X(K/\Q, C_1, C_2)) \gg \sum_{\chi \text{ symplectic}} \log A(\chi) \gg \log |d_K|.$$
\end{demo}

\begin{prop}\label{Esp1} Assume $\GRH$ and $\LIP$. If $K/\Q$ is tamely ramified then for any $3 \leq i \leq n$, $$\E(X(K/K_i, C_1^{(i)}, C_{-1}^{(i)})) = - 2^{n-1}(1-W_K) + 2^{i-1},$$ where $W_k = W_{k/\Q}$ is the root number of any symplectic character of $G$.
\end{prop}

\begin{demo} From Theorem \ref{VA} we have $$\E(X(K/K_i, C_1^{(i)}, C_{-1}^{(i)})) = \frac{|(C_{-1}^{(i)})^{1/2}|}{|C_{-1}^{(i)}|} - \frac{|(C_1^{(i)})^{1/2}|}{|C_1^{(i)}|} + 2 \sum_{\chi \neq \chi_0} (\chi(C_{-1}^{(i)})-\chi(C_1^{(i)})) \ord_{s=1/2} L(s, \chi, K/K_i).$$ Since we are assuming $\LIP$, Proposition \ref{Zero} shows that only symplectic characters, i.e. the $\psi_j^{(i)}$ with $1 \leq j \leq 2^{i-2}-1$ odd, contribute to the sum. By the same proposition, their root numbers are all equal to $W_K$ and the order of vanishing of their Artin $L$-functions is $0$ or $2^{n-i}$, so $$2 \sum_{\chi \neq \chi_0} (\chi(C_{-1}^{(i)})-\chi(C_1^{(i)})) \ord_{s=1/2} L(s, \chi, K/K_i) = \left\{\begin{array}{l}
0 \text{ if } W_K=1\\
2 \sum_{\chi \text{ symplectic}} (-4) \cdot 2^{n-i} = -2^n \text{ if } W_K = -1.
\end{array}\right.$$ In the proof of Lemma \ref{Orth}, we have observed that $-1$ has $2^{i-1}+2$ square roots in $\Gal(K/K_i)$, while $1$ has $2$. This yields $$\E(X(K/K_i, C_1^{(i)}, C_{-1}^{(i)})) = \left\{\begin{matrix}
2^{i-1} \text{ if } W_K=1\\
2^{i-1} - 2^{n} \text{ if } W_K = -1.
\end{matrix}\right.$$
\end{demo}

For the sake of completeness, we give below an exhaustive list of $\E(X(K/K_i, C_1, C_2))$ for each $3 \leq i \leq n$ and for all possible choices of distinct conjugacy classes $C_1, C_2$ of $G_i$.

\begin{prop}\label{EspQ} Assume $\GRH$ and $\LIP$. If $K/\Q$ is tamely ramified then for any $3 \leq i \leq n$,

$$\begin{array}{c|c|c|c}
C_1&C_2&\E(X(K/K_i, C_1, C_2))&\text{Conditions on the classes}\\
\hline
C_1^{(i)}&C_{-1}^{(i)}&- 2^{n-1}(1-W_K) + 2^{i-1}& none\\
\hline
C_1^{(i)}&C_{x_i^k}^{(i)}&- 2^{n-2}(1-W_K) + (-1)^k-1& 1 \leq k \leq 2^{i-2}-1\\
\hline
C_1^{(i)}&C_{y}^{(i)}/C_{x_iy}^{(i)}&-2^{n-2}(1-W_K)-2& none\\
\hline
C_{-1}^{(i)}&C_{x_i^k}^{(i)}&2^{n-2}(1-W_K) - 1+(-1)^{k} - 2^{i-1}& 1 \leq k \leq 2^{i-2}-1\\
\hline
C_{-1}^{(i)}&C_{y}^{(i)}/C_{x_iy}^{(i)}&2^{n-2}(1-W_K)-2-2^{i-1}& none\\
\hline
C_{x_i^k}^{(i)}&C_{x_i^l}^{(i)}&(-1)^l-(-1)^k& 1 \leq k, l \leq 2^{i-2}-1, k \neq l\\
\hline
C_{x_i^k}^{(i)}&C_{y}^{(i)}/C_{x_iy}^{(i)}&(-1)^{k+1}-1& 1 \leq k \leq 2^{i-2}-1\\
\hline
C_y^{(i)}&C_{x_iy}^{(i)}&0& none\\
\end{array}$$

\end{prop}

\begin{demo} The first row was computed in Proposition \ref{Esp1}. Each of the sums $$\sum_{\chi \neq \chi_0} \chi(C_{x_i^k}^{(i)}) \ord_{s=1/2} L(s, \chi, K/K_i)$$ for $1 \leq k \leq 2^{i-2}-1$ is zero, because Proposition \ref{Zero} shows they only involve symplectic characters, for which $\ord_{s=1/2} L(s, \chi, K/K_i)$ does not depend on $\chi$, and those sums reduce to $$\ord_{s=1/2} L(s, \psi_1^{(i)}, K/K_i) \sum_{\underset{j \text{ odd}}{j=1}}^{2^{i-2}-1} \left(\zeta_{i}^{jk} + \zeta_{i}^{-jk}\right),$$ which is zero by Lemma \ref{Root}. Therefore using Theorem \ref{VA}, we see that the remaining computations only involves counting square roots in $G_i \simeq \H_{2^i}$, and, when one of the conjugacy classes is $C_1^{(i)}$ or $C_{-1}^{(i)}$, the fact there are $2^{i-3}$ irreducible symplectic characters of $G_i$ with $L$-functions vanishing at $1/2$ to order $2^{n-i-1}(1-W_K)$. The last row is clear because the elements of $C_y^{(i)}$ and $C_{x_i y}^{(i)}$ have no square root in $G_i$ and symplectic characters vanish on those conjugacy classes. 
\end{demo}

The table of \ref{EspQ} that any Chebotarev race involving the classes $C_1^{(i)}$ or $C_{-1}^{(i)}$ is influenced by the root number of symplectic characters of $G$.

\subsection{Moments in dihedral Galois groups of order a power of two}

\label{Sect4.2}

We now turn to the case of number field extensions with dihedral Galois group of $2$-power order. The same methods as in the previous section are applied to get bounds on the moments of the random variables attached to the different Chebotarev races in such extensions. Those bounds will be used to prove Theorem \hyperlink{ThB}{B}.\\

Let $n \geq 2$ and let $L/\Q$ be a Galois extension with Galois group $D_{2^{n-1}} = \langle r,s \mid r^{2^{n-1}}=s^2=1, srs^{-1}=r^{-1} \rangle.$ Let $\Q = L_n \subset L_{n-1} \subset \dots \subset L_3 \subset L$ be the subextensions such that $G_i := \langle r^{2^{n-i}}, s \rangle$ and $L_i=L^{G_i}$. We will denote $r^{2^{n-i}}$ by $r_i$. Conjugacy classes in $G_i$ will be denoted, following our notations in Lemma \ref{ConjD}, by $C_1^{(i)}, C_{-1}^{(i)}, C_{r_i^k}^{(i)}, C_s^{(i)}$ and $C_{r_i s}^{(i)}$. The estimates from the previous section are proved similarly in the dihedral case, so we state them without proof. Since each irreducible character of $D_{2^{n-1}}$ is orthogonal by Lemma \ref{Orth}, and since we will be working under the hypothesis $\LI$, no considerations on root numbers are involved.

\begin{prop}\label{VarD} Assume $\GRH$ and $\LIM$. Then for any $2 \leq i \leq n$ and any distinct conjugacy classes $C_1$, $C_2$ of $G_i$ such that $\{C_1, C_2\} \neq \{C^{(i)}_{s}, C^{(i)}_{r_i s}\}$,, $$\Var(X(L/L_i, C_1, C_2)) \ll \log |d_L|.$$ 
Moreover if $L/\Q$ is tamely ramified and $\{C_1, C_2\} \not \subset \{C^{(i)}_{r_i^{2^{i-3}}}, C^{(i)}_{s}, C^{(i)}_{r_i s}\}$ then for any $2 \leq i \leq n,$ $$\Var(X(L/L_i, C_1, C_2)) \gg \frac{\log |d_L|}{16^{n}},$$ and if $\{C_1, C_2\} \cap \{C_{r_i^k}^{(i)} \mid 1 \leq k \leq 2^{i-2}-1\} = \emptyset$ this can be improved to $$\Var(X(L/L_i, C_1, C_2)) \gg \log |d_L|.$$
\end{prop}

\begin{prop}\label{EspD} Assume $\GRH$ and $\LI$. For any $3 \leq i \leq n$,

$$\begin{array}{c|c|c|c}
C_1&C_2&\E(X(L/L_i, C_1, C_2))&\text{Conditions on the classes}\\
\hline
C_1^{(i)}&C_{-1}^{(i)}&-2^{i-1}& none\\
\hline
C_1^{(i)}&C_{r_i^k}^{(i)}&- 2^{i-1} + (-1)^k - 1& 1 \leq k \leq 2^{i-2}-1\\
\hline
C_1^{(i)}&C_s^{(i)}/C_{r_is}^{(i)}&-2^{i-1} - 2& none\\
\hline
C_{-1}^{(i)}&C_{r_i^k}^{(i)}&(-1)^k-1& 1 \leq k \leq 2^{i-2}-1\\
\hline
C_{-1}^{(i)}&C_{s}^{(i)}/C_{r_is}^{(i)}&-2& none\\
\hline
C_{r_i^k}^{(i)}&C_{r_i^l}^{(i)}&(-1)^l-(-1)^k& 1 \leq k, l \leq 2^{i-2}-1, k \neq l\\
\hline
C_{r_i^k}^{(i)}&C_{s}^{(i)}/C_{r_is}^{(i)}&(-1)^{k+1}-1& 1 \leq k \leq 2^{i-2}-1\\
\hline
C_s^{(i)}&C_{r_is}^{(i)}&0& none\\
\end{array}$$
\end{prop}

\subsection{Construction of the towers}

\label{Sect4.3}

We now construct towers of dihedral and generalized quaternion extensions of $\Q$ which are tamely ramified, and with controlled discriminants, properties which will enable us to apply effectively the results from sections 4.1 and 4.2.

\begin{prop}\label{ConsD} Assume $\GRH$ for the Dedekind zeta functions of the number fields\\ $\Q(\sqrt[5 \cdot 2^{n-1}]{\varepsilon}, \mu_{5 \cdot 2^{n-1}})$, where $\mu_k$ denotes a primitive $k$-th root of unity and $\varepsilon = \frac{3+\sqrt{5}}{2}$. There exists a sequence $(\mathcal{D}_n)_{n \geq 3}$ of number fields such that for any $n \geq 3$ :
\begin{itemize}
\item[i)] The extension $\mathcal{D}_n/\Q$ is tamely ramified and Galois with Galois group isomorphic to the dihedral group $D_{2^{n-1}}$ of order $2^n$.
\item[ii)] $\Q(\sqrt{5}) \subset \mathcal{D}_n$.
\item[iii)] $2^n \ll \log |d_{\mathcal{D}_n/\Q}| \ll n2^n$.
\end{itemize}
\end{prop}

\begin{demo} Fix $K = \Q(\sqrt{5})$, let $n \geq 3$ and let $\varepsilon = \frac{3+\sqrt{5}}{2}$ be the fundamental totally positive unit of $K$. Using class field theory, it is shown in the proof of \cite[Theorem 3.2]{Pla} that there exists a number field $\mathcal{D}_n$ which has Galois group over $\Q$ isomorphic to $D_{2^{n-1}}$, containing $K$ and such that only $5$ and another odd prime $p$ ramify in $\mathcal{D}_n$. The prime number $p$ is chosen so that $p$ splits in the extension $M_n := \Q(\sqrt[5 \cdot 2^{n-1}]{\varepsilon}, \mu_{5 \cdot 2^{n-1}})$.

Using the bound for the least prime ideal in the Chebotarev density theorem stated in \cite[(1.2)]{LMO}, we can find such a $p$ satisfying $$p \ll (\log |d_{M_n}|)^2.$$ Now, $M_n$ is obtained in at most $2n$ steps of adjoining square roots of algebraic units, starting from the field $\Q(\sqrt[5]{\varepsilon}, \mu_5)$. If $M/L$ is one of those steps, we have $M=L(\sqrt{\eta})$, where $\eta$ is a unit of $L$, and the relative discriminant $D_{M/L}$ has to divide the discriminant of $X^2 - \eta$, which is $4\eta\mathcal{O}_L = 4\mathcal{O}_L$. Since $$|d_M| = N_{L/\Q}(D_{M/L}) |d_L|^{[M:L]} \leq 4^{[L:\Q]} |d_L|^2,$$ we find $$|d_M|^{\frac{1}{[M:\Q]}} \leq 2 |d_L|^{\frac{1}{[L:\Q]}}.$$ Using this iteratively on the (at most) $2n$ steps, we find that $$|d_{M_n}|^{\frac{1}{[M_n : \Q]}} \ll 2^{2n}.$$ Using the same reasoning, we see that $$[M_n : \Q] \ll 4^n,$$ so finally we can choose the prime $p$ so that $$p \ll (n[M_n:\Q])^2 \ll n^2 16^n.$$

As in the proof of Proposition \ref{VarQ}, $$|d_{\mathcal{D}_n}| = 5^{(e_5-1)f_5g_5} \times p^{(e_p-1)f_pg_p} \leq (5p)^{2^n},$$ where $e_q$ (respectively $f_q$, $g_q$) denotes the ramification index of the prime number $q$ (respectively the residual degree of $q$, the number of primes above $q$). This shows that $$\log |d_{\mathcal{D}_n}| \ll 2^n \log p \ll n2^n.$$ The lower bound on the discriminant simply comes from Minkowski's bound (\cite[P.120]{Lan}) $$|d_{\mathcal{D}_n}| \geq \frac{2^{n2^n}}{(2^n)!} \left(\frac{\pi}{4}\right)^{2^{n-1}}.$$ Since $[\mathcal{D}_n : \Q] = 2^n$ and $d_{\mathcal{D}_n}$ is odd, $\mathcal{D}_n/\Q$ is tamely ramified.
\end{demo}

\begin{rmq} The quadratic field $\Q(\sqrt{5})$ plays no particular role in our construction. For our purpose, we could replace the integer $5$ by any positive square-free integer $d$ which is $1$ mod $4$.
\end{rmq}

\begin{prop}\label{ConsQ} Assume $\GRH$. There exists two sequences $(\mathcal Q_n^+)_{n \geq 3}$ and $(\mathcal Q_n^-)_{n \geq 3}$ of number fields such that for any $n \geq 3$ :
\begin{itemize}
\item[i)] The extension $\mathcal Q_n^{\pm}/\Q$ is tamely ramified and Galois with Galois group isomorphic to the generalized quaternion group $\H_{2^n}$.
\item[ii)] $\Q(\sqrt{5}) \subset \mathcal Q^{\pm}_n$.
\item[iii)] $2^n \ll \log |d_{\mathcal Q_n^{\pm}/\Q}| \ll n2^n$.
\item[iv)] $W_{\mathcal Q_n^+} = 1$ and $W_{\mathcal Q_n^-}=-1$.
\end{itemize}
\end{prop}

\begin{demo} We use the same notations as in the proof of Proposition \ref{ConsD} : we let $K = \Q(\sqrt{5})$ and $\varepsilon = \frac{3+\sqrt{5}}{2}$ be the fundamental totally positive unit of $K$. We use again a theorem of Fröhlich \cite[Chapter V, Proposition 3.1]{Fro3}. It states that for $e \in \{-1, 1\}$, the set of prime numbers $p$ such that there exists a number field $N[p]$ such that the only ramified primes in $N[p]$ are $5$ and $p$, satisfying $\Gal(N[p]/\Q) \simeq \H_8$, $\Q(\sqrt{5}) \subset N[p]$ and $W_{N[p]} = e$ have positive density. It is shown by translating the last two conditions into an arithmetic condition on $p$, and then using Chebotarev's density theorem. The arithmetic condition amounts to prescribing the Frobenius conjugacy class of $p$ in $\Gal(M_n'/\Q)$, where $M_n' = K(\mu_{2^{n-1}}, \sqrt[2^{n-2}]{\varepsilon})$, with $\mu_{2^{n-1}}$ a primitive $2^{n-1}$-th root of unity. As in the proof of Proposition \ref{ConsD}, we apply, under $\GRH$, the bound on the least prime ideal in Chebotarev's density theorem of \cite{LMO} to choose $$p \ll (\log |d_{M_n'}|)^2 \ll n^216^n.$$ Call $\mathcal{Q}_n^{\pm}$ the number field constructed this way (the superscript $\pm$ indicating which root number was prescribed). Since only $5$ and $p$ ramify in $\mathcal Q_n^{\pm}$, we get as in the proof of Proposition \ref{ConsD} $$\log |d_{\mathcal Q_n^{\pm}/\Q}| \ll n2^n.$$ The lower bound on the discriminant follows once again from Minkowski's bound, and tame ramification follows from the fact that $d_{\mathcal Q_n^{\pm}}$ is odd.
\end{demo}

\begin{rmq} If we do not want to specify the root numbers in our quaternion extensions of $\Q$, we could use the methods of \cite{DaMa} to construct extensions of $\Q$ satisfying i), ii) and iii) using the extensions $\mathcal{D}_n/\Q$ of Proposition \ref{ConsD}. Indeed, it is shown in \cite{DaMa} that, if $Q/\Q$ is a tamely ramified extension with Galois group $\H_8$ and $\Q(\sqrt{5}) \subset Q$ (again, the number $5$ has no particular significance), then the composite field $Q\mathcal{D}_n$ contains a subfield $Q_n$ with Galois group $\H_{2^n}$ over $\Q$, and the upper bound on $\log |d_{\mathcal{D}_n}|$ then implies a similar bound on $\log |d_{Q_n}|$. Tame ramification follows from the tame ramification of $Q/\Q$ and $\mathcal{D}_n/\Q$, and the lower bound on the discriminant from Minkowski's bound.
\end{rmq}

\subsection{A large deviation result for Chebyshev's bias}

\label{Sect4.4}

In this section we prove a lower bound in the context of Theorem \ref{DevFJ}. We will make use of the following bounds on sums of zeros of Artin $L$-functions. The first one is the so-called Riemann-Von Mangoldt formula, stated in \cite{IwKo}, in the particular case of Artin $L$-functions, for which the conductor of $s \mapsto L(s, \chi, L/K)$ is simply $A(\chi)$ and its analytic conductor is bounded by $A(\chi)(|s|+4)^{[K:\Q]\chi(1)}$ (see \cite[§5.13]{IwKo}).

\begin{lem} \label{N} Let $\chi$ be an irreducible character of $G$, and for any $T > 0$, define $N(T, \chi)$ to be the number of zeros $\rho = \beta + i \gamma$ of $s \mapsto L(s, \chi, L/K)$ with $0 \leq \beta \leq 1$ and $0 < \gamma \leq T$. Then we have $$N(T, \chi) = \frac{T}{2\pi} \log\left(A(\chi) \left(\frac{T}{2 \pi e}\right)^{[K:\Q] \chi(1)}\right) + O\left(\log\left(A(\chi) (T+4)^{[K:\Q] \chi(1)}\right)\right).$$
\end{lem}

We note that, in the above Lemma, we are counting zeros $\beta + i \gamma$ with $0 < \gamma \leq T$ and not $|\gamma| \leq T$, so the main term is half the one in \cite[Theorem 5.8]{IwKo}, while the number of zeros is taken into account in the error term.

\begin{lem}[\cite{FiJo}, Lemma 5.4]\label{Som} Assume Artin's conjecture. Then for any irreducible character $\chi$ of $G$ and for $T \geq 1$, we have $$\sum_{0 < \gamma_{\chi} \leq T} \frac{1}{\sqrt{\frac{1}{4} + \gamma_{\chi}^2}} =  \frac{\log T}{2\pi} \log\left(A(\chi) \left(\frac{T^{1/2}}{2 \pi e}\right)^{[K:\Q] \chi(1)}\right) + O\left(\log\left(A(\chi) (T+4)^{[K:\Q] \chi(1)}\right)\right).$$
\end{lem}

In order to state our improvement of Theorem \ref{DevFJ}, we first recall the large deviation result of Montgomery and Odlyzko which was used to derive Theorem \ref{DevFJ}, and will be used to prove Theorem \ref{Dev}.

\begin{thm}[\cite{MoOd}, Theorem 2]\label{DevMO} Let $(W_n)_{n \geq 1}$ be a family of independent real random variables such that \begin{itemize}
\item[i)] For all $n \geq 1$, $\E(W_n)=0$.
\item[ii)] For all $n \geq 1$, $|W_n| \leq 1$ a.s.
\item[iii)] There exists $c > 0$ such that for all $n \geq 1$, $\E(W_n^2) > c$.
\end{itemize}
Let $(r_n)_n$ be a real sequence decreasing to zero such that $\displaystyle \sum_{n \geq 1} r_n^2 < +\infty$ and let $W = \displaystyle \sum_{n \geq 1} r_n W_n$. Let $V \geq 0$ and $\alpha > 0$. If $\displaystyle \sum_{r_n \geq \alpha} r_n \leq \frac{V}{2}$ then $$\P(W \geq V) \leq \exp\left(-\frac{1}{16} V^2 \left(\sum_{r_n < \alpha} r_n^2\right)^{-1}\right),$$
If $\displaystyle \sum_{r_n \geq \alpha} r_n \geq 2V$ then $$\P(W \geq V) \geq a_1 \exp\left(-a_2 V^2 \left(\sum_{r_n < \alpha} r_n^2\right)^{-1}\right),$$ where $a_1$ and $a_2$ are positive constants depending only on $c$.
\end{thm}

The following theorem provides a lower bound for $1 - \delta$, which was not obtained in \cite{FiJo}, in the context of Theorem \ref{DevFJ}. The proof proceeds by taking into account distinct characters of $G$ whose corresponding induced characters to $G^+$ are not orthogonal, which in general leads to complications in estimating Chebyshev's bias. This is because, in Theorem \ref{VA}, the means are expressed using values of characters of $G$, while the variances involve values of characters of $G^+$.\\

Before stating the result, we set a few notations for readability. If $\chi \in \Irr(G)$ and $\lambda \in \Irr(G^+)$ we will write $\lambda \mid \chi$ if $\langle \lambda, \Ind_G^{G^+} \chi \rangle \neq 0$, that is if $\lambda$ is a component of the character $\Ind_G^{G^+} \chi$ of $G^+$, and $\lambda \nmid \chi$ otherwise. For $\chi, \chi' \in \Irr(G)$, we write $(\chi, \chi')=1$ when $\langle \Ind_G^{G^+} \chi, \Ind_G^{G^+} \chi'\rangle = 0$, \textit{i.e.} when $\lambda \mid \chi$ implies $\lambda \nmid \chi'$.

\begin{thm}\label{Dev} Assume $\GRH, \LI$ and Artin's Conjecture. Let $$S = \{\chi \in \Irr(G) \mid \forall \chi' \in \Irr(G) \setminus \{\chi\}, (\chi, \chi')=1\}$$ and $R = \Irr(G) \setminus (S \cup \{\chi_0\})$. Let $b_1 := \max_{\chi \in R} \chi(1)$, $b_2 := |R|$ and $M=\max_{\chi \in R} \ord_{s=1/2} L(s, \chi, L/K)+1$. There exist absolute constants $c_1, c_2 > 0$ such that for any conjugacy classes $C_1, C_2$ of $G$ satisfying $C_1^+ \neq C_2^+$ and $B(L/K, C_1, C_2) > 0$, we have $$c_1 \exp\left(-c_2 Q(C_1, C_2) B(L/K, C_1, C_2)^2\right) \leq 1 - \delta(L/K, C_1, C_2)$$ where $Q(C_1, C_2) = \max\left(e^{C\sqrt{\frac{Mb_1 b_2}{\lambda^*(1)b_3(C_1, C_2)}}}, C\frac{b_3(C_1, C_2)}{b_4(C_1, C_2)}, C\right)$ for some absolute constant $C > 0$ and where $\lambda^* \in \Irr(G^+)$ is such that $|\lambda^*(C_2^+)-\lambda^*(C_1^+)| = \max_{\lambda \in \Irr(G^+)} |\lambda(C_2^+)-\lambda(C_1^+)| =: b_3(C_1, C_2) > 0$ and $b_4(C_1, C_2) = \min_{\lambda \in \Irr(G^+), \lambda(C_2^+) \neq \lambda(C_1^+)} |\lambda(C_2^+)-\lambda(C_1^+)| > 0$.
\end{thm}

\begin{demo} We introduce the random variable $W = X(L/K, C_1, C_2) - \E(X(L/K, C_1, C_2))$. It obviously satisfies the hypotheses of Theorem \ref{DevMO} with $(r_n)_n$ the ordered sequence of non-zero $\frac{2|\lambda(C_2^+) - \lambda(C_1^+)|}{\sqrt{\frac{1}{4} + \gamma_{\lambda}^2}}$, $\lambda$ varying in $\Irr(G^+)$ and $\gamma_{\lambda}$ in $\Gamma_{L/\Q, \lambda}$, and $(W_n)_n$ the accordingly ordered sequence of $X_{\gamma}$. Indeed, the series $$\sum_{\gamma_{\lambda} > 0} \frac{1}{\frac{1}{4} + \gamma_{\lambda}^2}$$ converges for any $\lambda \in \Irr(G^+)$.\\

We set $V = \E(X(L/K, C_1, C_2))$ and look for $\alpha > 0$ such that $\displaystyle \sum_{r_n \geq \alpha} r_n \geq 2V$. By Theorem \ref{VA}, we have \begin{align*}\E(X(L/K, C_1, C_2)) &= \frac{|C_2^{1/2}|}{|C_2|} - \frac{|C_1^{1/2}|}{|C_1|} + z(C_2) - z(C_1)\\ &= \sum_{\chi \in \Irr(G)} (\chi(C_2)-\chi(C_1))(\varepsilon_2(\chi) + 2 \ord_{s=1/2} L(s, \chi, L/K)).\end{align*} If $\chi \in S$, then for some $\lambda \in \Irr(G^+)$ we have, by a quick computation using the Frobenius reciprocity formula, \begin{align*}|\lambda(C_2^+) - \lambda(C_1^+)| = |\lambda_{\mid G}(C_2) - \lambda_{\mid G}(C_1)| &= \left|\sum_{\chi' \in \Irr(G)} (\chi'(C_2) - \chi'(C_1)) \langle \lambda_{\mid G}, \chi'\rangle\right|\\ &= |\chi(C_2)-\chi(C_1)| \langle \lambda, \Ind_G^{G^+} \chi\rangle.\end{align*} From the factorisation $$L(s, \chi, L/K) = \prod_{\lambda \in \Irr(G^+)} L(s, \lambda, L/\Q)^{\langle \lambda, \Ind_G^{G^+} \chi \rangle},$$ we see that $$\ord_{s=1/2} L(s, \chi, L/K) = \sum_{\lambda \in \Irr(G^+)} \langle \lambda, \Ind_G^{G^+} \chi \rangle \ord_{s=1/2} L(s, \lambda, L/\Q) \leq M_0 \sum_{\lambda \in \Irr(G^+)} \langle \lambda, \Ind_G^{G^+} \chi \rangle$$ (recall that $M_0$ comes from the hypothesis $\LI$). Therefore, \begin{align*}\sum_{\chi \in S} (\chi(C_2)-\chi(C_1))(\varepsilon_2(\chi) + 2 \ord_{s=1/2} L(s, \chi, L/K)) &\leq \sum_{\chi \in S} |\chi(C_2)-\chi(C_1)| (1+2M_0\sum_{\lambda \mid \chi} \langle \lambda, \Ind_G^{G^+} \chi \rangle)\\ &=(1+2M_0) \sum_{\lambda \in S^+} |\lambda(C_2^+)-\lambda(C_1^+)|\end{align*} where $S^+ = \{\lambda \in \Irr(G^+) \mid \exists \chi \in S, \lambda \mid \chi\}$.  We now see that it is enough to have $$\sum_{0 < \gamma_{\lambda} \leq T_0(\lambda)} \frac{1}{\sqrt{\frac{1}{4} + \gamma_{\lambda}^2}} \geq (2+4M_0)$$ for any $\lambda \in S^+$ and some $T_0(\lambda) > 0$, to bound from above the contribution of characters of $S$ in $2V$. It remains to bound from above the sum $$2\sum_{\chi \in R} (\chi(C_2)-\chi(C_1))(\varepsilon(\chi) + 2\ord_{s=1/2} L(s, \chi, L/K)) \leq 4M \sum_{\chi \in R} |\chi(C_2)-\chi(C_1)|.$$ By definition of $b_1$ and $b_2$, this last sum is $\leq 8M b_1 b_2$, so we choose $T_0(\lambda^*) > 0$ such that $$\sum_{0 < \gamma_{\lambda^*} \leq T_0(\lambda^*)} \frac{1}{\sqrt{\frac{1}{4} + \gamma_{\lambda^*}^2}} \geq \frac{8M b_1 b_2}{b_3(C_1, C_2)}.$$

To do so, we use Lemma \ref{Som} (for the extension $L/\Q$) to get $$\sum_{0 < \gamma_{\lambda} \leq T} \frac{1}{\sqrt{\frac{1}{4} + \gamma_{\lambda}^2}} = \frac{\log T}{2\pi} \log\left(A(\lambda) \left(\frac{T^{1/2}}{2 \pi e}\right)^{\lambda(1)}\right) + O\left(\log\left(A(\lambda) (T+4)^{\lambda(1)}\right)\right)$$ for any $\lambda \in \Irr(G^+)$. A simple computation shows that it is enough to choose $T_0(\lambda^*) = \max\left((2\pi e)^4, e^{\sqrt{\frac{64\pi Mb_1 b_2}{\lambda^*(1)b_3(C_1, C_2)}}}, e^{16 \pi C'}\right)$ to ensure $$\sum_{0 < \gamma_{\lambda^*} < T_0(\lambda^*)} \frac{1}{\sqrt{\frac{1}{4} + \gamma_{\lambda^*}^2}} \geq \frac{8Mb_1 b_2}{b_3(C_1, C_2)},$$ where $C' > 0$ is such that $$\left|\sum_{0 < \gamma_{\lambda^*} \leq T} \frac{1}{\sqrt{\frac{1}{4} + \gamma_{\lambda^*}^2}} - \frac{\log T}{2\pi} \log\left(A(\lambda^*) \left(\frac{T^{1/2}}{2 \pi e}\right)^{\lambda^*(1)}\right)\right| \leq C' \log\left(A(\lambda^*) (T+4)^{\lambda^*(1)}\right),$$ while we have $$\sum_{0 < \gamma_{\lambda} \leq A} \frac{1}{\sqrt{\frac{1}{4} + \gamma_{\lambda}^2}} \geq (2+4M_0)$$ for any $\lambda \in \Irr(G^+)$ and some absolute constant $A > 0$, so we set $T_0(\lambda)=A$ for every $\lambda \neq \lambda^*$.\\

We now turn to the choice of $\alpha \geq 0$ such that $$\sum_{r_n \geq \alpha} r_n \geq 2V.$$ It is enough to have $$\sum_{r_n \geq \alpha} r_n \geq \sum_{\lambda \in \Irr(G^+)} \sum_{0 < \gamma_{\lambda} \leq T_0(\lambda)} \frac{2|\lambda(C_2^+) - \lambda(C_1^+)|}{\sqrt{\frac{1}{4} + \gamma_{\lambda}^2}},$$ so we require that $$0 < \gamma_{\lambda} \leq T_0(\lambda) \Rightarrow \gamma_{\lambda} \leq \sqrt{\frac{4|\lambda(C_2^+) - \lambda(C_1^+)|^2}{\alpha^2} - \frac{1}{4}} =: R_{\alpha, \lambda}.$$ We choose $\alpha$ to be the minimal (positive) value of $\frac{2|\lambda(C_2^+) - \lambda(C_1^+)|}{\sqrt{\frac{1}{4} + T_0(\lambda)^2}}$, attained for some $\lambda_m \in \Irr(G^+)$ (recall that we discarded the characters $\lambda \in \Irr(G^+)$ such that $\lambda(C_2^+) = \lambda(C_1^+)$ in our definition of the $r_n$'s). For this choice of $\alpha$ we therefore have $\sum_{r_n \geq \alpha} r_n \geq 2V$.\\

Theorem \ref{DevMO} now yields $$\P(W \geq V) \geq a_1 \exp\left(-a_2 V^2 \left(\sum_{r_n < \alpha} r_n^2\right)^{-1}\right).$$ We have \begin{align*}
\sum_{r_n < \alpha} r_n^2 &= \sum_{\lambda \in \Irr(G^+)} \sum_{\gamma_{\lambda} > R_{\alpha, \lambda}} \frac{4|\lambda(C_2^+) - \lambda(C_1^+)|^2}{\frac{1}{4} + \gamma_{\lambda}^2}\\ &\geq \sum_{\lambda \in \Irr(G^+)} \sum_{2 R_{\alpha, \lambda} \geq \gamma_{\lambda} > R_{\alpha, \lambda}} \frac{4|\lambda(C_2^+) - \lambda(C_1^+)|^2}{\frac{1}{4} + \gamma_{\lambda}^2}\\ &\geq \sum_{\lambda \in \Irr(G^+)} \frac{4|\lambda(C_2^+) - \lambda(C_1^+)|^2}{\frac{1}{4} + 4R_{\alpha, \lambda}^2}(N(2R_{\alpha, \lambda}, \lambda) - N(R_{\alpha, \lambda}, \lambda))
\end{align*}

As in the proof of \cite[Lemma 4.3]{FiJo} we see, using Lemma \ref{N}, that $$N(2R_{\alpha, \lambda}, \lambda) - N(R_{\alpha, \lambda}, \lambda) \gg R_{\alpha, \lambda} \log A(\lambda).$$ Thus, $$\sum_{r_n < \alpha} r_n^2 \gg \frac{1}{\max_{\lambda \in \Irr(G^+)} R_{\alpha, \lambda}} \sum_{\lambda \in \Irr(G^+)} |\lambda(C_2^+) - \lambda(C_1^+)|^2 \log A(\lambda).$$ Combining Theorem \ref{VA} and Lemma \ref{B0}, we finally get $$\sum_{r_n < \alpha} r_n^2 \gg\frac{1}{\max_{\lambda \in \Irr(G^+)} R_{\alpha, \lambda}} \Var(X(L/K, C_1, C_2)).$$ We can now conclude that $$\P(W \geq V) \geq a_1 \exp\left(-a_3 (\max_{\lambda \in \Irr(G^+)} R_{\alpha, \lambda}) B(L/K, C_1, C_2)^2 \right)$$ for some absolute constant $a_3 > 0$. But, recalling the choice we made for $\alpha$, $$\max_{\lambda \in \Irr(G^+)} R_{\alpha, \lambda} = \sqrt{\left(\frac{1}{4} + T_0(\lambda_m)^2\right) \frac{|\lambda^*(C_2^+) - \lambda^*(C_1^+)|^2}{|\lambda_m(C_2^+) - \lambda_m(C_1^+)|^2}}.$$ To determine the size of this quantity, we consider two cases : either $\lambda_m=\lambda^*$, in which case we find $\max_{\lambda \in \Irr(G^+)} R_{\alpha, \lambda} \asymp T_0(\lambda^*) \asymp e^{C \sqrt{\frac{Mb_1 b_2}{\lambda^*(1)b_3(C_1, C_2)}}}$, or $\lambda_m \neq \lambda^*$, in which case $\max_{\lambda \in \Irr(G^+)} R_{\alpha, \lambda} \asymp \frac{b_3(C_1, C_2)}{b_4(C_1, C_2)}$ since $T_0(\lambda_m)$ was chosen as an absolute constant.\\

It simply remains to note that, by symmetry of $W$ about $0$, we have $$\P(W \geq V) = \P(W \leq -V) = 1 - \P(X(L/K, C_1, C_2) > 0) = 1 - \delta(L/K, C_1, C_2).$$
\end{demo}

Note that when $K=\Q$ in the previous theorem, we have $R = \emptyset$ so that $Q(C_1, C_2)$ can be taken to be $C\left(\frac{b_3(C_1, C_2)}{b_4(C_1, C_2)}+1\right),$ where $C > 0$ is an absolute constant. In our applications (Theorems \ref{TabD} and \ref{TabQ}), the quantities $b_1, b_2, M$ and $\frac{b_3(C_1, C_2)}{b_4(C_1, C_2)}$ will be bounded.

\subsection{Estimates on the bias in the towers}

\label{Sect4.5}

Using the number fields constructed in Proposition \ref{ConsD} and Proposition \ref{ConsQ}, we can now state the following result, a more exhaustive version of Theorem \hyperlink{ThB}{B}.

\begin{thm}\label{TabD} Assume $\GRH$ and $\LI$ for the number fields $\mathcal{D}_n$ of Proposition \ref{ConsD}. There exist absolute constants $c_1, c_2, c_3 > 0$ such that for any $n \geq 3$, the following hold :

$$\begin{array}{c|c|c|c}
C_a & C_b & \text{Estimate on } \delta(\mathcal{D}_n/\Q, C_1, C_2) & \text{Condition on the classes}\\
\hline
C_1 & C_{-1}, C_s, C_{rs} & c_1\exp(-c_22^n) < \delta(\mathcal{D}_n/\Q, C_1, C_{b}) < \exp\left(-c_3 \frac{2^n}{n}\right) & none\\
\hline
C_1 & C_{r^k} & c_1\exp(-c_3 32^n) < \delta(\mathcal{D}_n/\Q, C_1, C_{r^k}) < \exp\left(-c_4 \frac{2^n}{n}\right) & none\\
\hline
C_{-1} & C_{r^k} & \frac{1}{2} & k \text{ even}\\
\hline
C_{-1} & C_s, C_{rs} & 0 < \frac{1}{2} - \delta(\mathcal{D}_n/\Q, C_{-1}, C_b) \ll \frac{1}{2^n} & none\\
\hline
C_{r^k} & C_{r^l} & \frac{1}{2} & k = l \text{ mod } 2\\
\hline
C_{r^k} & C_s, C_{rs} & \frac{1}{2} & k\text{ odd}\\
\hline
C_s & C_{rs} & \frac{1}{2} & none
\end{array}$$
\end{thm}

\begin{demo} For the first part of the theorem, the bounds obtained in Proposition \ref{VarD} (for the variance), Proposition \ref{EspD} (for the mean) and Proposition \ref{ConsD} (for the discriminant) show that $$\frac{2^n}{n} \ll B(\mathcal{D}_n/\Q, C_1, C_{b})^2 = \frac{\E(X(\mathcal{D}_n/\Q, C_1, C_{b}))^2}{\Var(X(\mathcal{D}_n/\Q, C_1, C_{b}))} \ll 2^n.$$ We combine this estimate with Theorems \ref{DevFJ} and \ref{Dev} (and the fact $\E(X(\mathcal{D}_n/\Q, C_1, C_{b})) < 0$) to get $$c_1 \exp(-c_2 2^n) < \delta(\mathcal{D}_n/\Q, C_1, C_{-1}) < \exp\left(-c_3 \frac{2^n}{n}\right),$$ for some absolute constants $c_1, c_2, c_3 > 0$ (which may differ from the constants of Theorems \ref{DevFJ} and \ref{Dev}, but by an absolute factor).

The proof is similar for $\delta(\mathcal{D}_n/\Q, C_1, C_{r^k})$, the only difference being the lower bound $$\Var(X(\mathcal{D}_n/\Q, C_1, C_{r^k})) \gg \frac{\log |d_{\mathcal{D}_n}|}{16^n} \gg \frac{1}{8^n}$$ from Proposition \ref{VarD} and \ref{ConsD}. To deal with $\delta(\mathcal{D}_n/\Q, C_{-1}, C_b)$, $b=s$ or $b=rs$, we simply apply Theorem \ref{TCL}, together with Proposition \ref{VarD}.

Finally, each case in which $\delta(C_1, C_2, L/\Q) = \frac{1}{2}$ comes from the fact that $\E(X(C_1, C_2, L/\Q))=0$.
\end{demo}

\begin{rmq} 
We could not produce estimates for $\delta(\mathcal{D}_n/\Q, C_{r^k}, C_b)$ for $b=-1$ (when $k$ is odd), $b=r^l$ (when $k$ and $l$ do not have the same parity) and $b=s$ or $b=rs$ (when $k$ is even) because we do not have good enough bounds on $\Var(X(\mathcal{D}_n/\Q, C_{r^k}, C_b))$ to conclude that $B(\mathcal{D}_n/\Q, C_{r^k}, C_b)$ is large or not.\\
\end{rmq}

We now turn our attention to the extensions $\mathcal Q_n^{\pm}/\Q$ built in Proposition \ref{ConsQ}. The proof of the next theorem is the same as for Theorem \ref{TabD}, except that the value of $W_{\mathcal Q_n^{\pm}}$ determines in some cases the class towards which there is a bias. This is a more exhaustive version of Theorem \hyperlink{ThC}{C}.

\begin{thm}\label{TabQ} Assume $\GRH$ and $\LIP$ for the number fields $\mathcal{Q}_n^{\pm}$ of Proposition \ref{ConsQ}. There exist absolute constants $c_1, c_2, c_3 > 0$ such that for any $n \geq 3$, denoting $\mathcal{Q}_n^{\pm}$ by $\mathcal{Q}_n$, the following hold :

$$\begin{array}{c|c|c|c}
C_a & C_b & \text{Estimate on } \delta(\mathcal{Q}_n/\Q, C_a, C_b) & \text{Conditions}\\
\hline
C_{-1} & C_{1} & c_1\exp(-c_22^n) < \left|\frac{1-W_{\mathcal Q_n}}{2} - \delta(\mathcal Q_n/\Q, C_{-1}, C_{b})\right| < \exp\left(-c_3 \frac{2^n}{n}\right) & none\\
\hline
C_{-1} & C_y, C_{xy} & c_1\exp(-c_22^n) < \delta(\mathcal Q_n/\Q, C_{-1}, C_{b}) < \exp\left(-c_3 \frac{2^n}{n}\right) & W_{\mathcal{Q}_n}=1\\
\hline
C_{-1} & C_y, C_{xy} & 0 < \frac{1}{2} - \delta(\mathcal Q_n/\Q, C_{-1}, C_{b}) \ll \frac{1}{2^{n/3}} & W_{\mathcal{Q}_n}=-1\\
\hline
C_{-1} & C_{x^k} & c_1\exp(-c_2 32^n) < \delta(\mathcal Q_n/\Q, C_{-1}, C_{x^k}) < \exp\left(-c_3 \frac{2^n}{n}\right) & W_{\mathcal{Q}_n}=1\\
\hline
C_{-1} & C_{x^k} & \frac{1}{2} & W_{\mathcal{Q}_n}=-1 \text{ and } k \text{ even}\\
\hline
C_{1} & C_{x^k} & \frac{1}{2} & W_{\mathcal Q_n}=1 \text{ and } k \text{ even}\\
\hline
C_{1} & C_{x^k} & c_1\exp(-c_2 32^n) < \delta(\mathcal Q_n/\Q, C_{1}, C_{x^k}) < \exp\left(-c_3 \frac{2^n}{n}\right) & W_{\mathcal Q_n} = -1\\
\hline
C_{1} & C_{y}, C_{xy} & 0 < \frac{1}{2} - \delta(\mathcal Q_n/\Q, C_1, C_b) \ll \frac{1}{2^{n/3}} & W_{\mathcal Q_n}=1\\
\hline
C_{1} & C_y, C_{xy} & c_1\exp(-c_2 2^n) < \delta(\mathcal Q_n/\Q, C_{1}, C_{b}) < \exp\left(-c_3 \frac{2^n}{n}\right) & W_{\mathcal Q_n} = -1\\
\hline
C_{x^k} & C_{x^l} & \frac{1}{2} & k = l \text{ mod } 2\\
\hline
C_{x^k} & C_y, C_{xy} & \frac{1}{2} & k\text{ odd}\\
\hline
C_y & C_{xy} & \frac{1}{2} & none
\end{array}$$
\end{thm}

\begin{rmq} As in Theorem \ref{TabD}, we could not produce bounds for $\delta(\mathcal Q_n/\Q, C_{x^k}, C_b)$ for $b=1$ (when $k$ is odd and $W_{\mathcal Q_n}=1$), $b=-1$ (when $k$ is odd and $W_{\mathcal{Q}_n}=-1$), $b=x^l$ (when $k$ and $l$ do not have the same parity) and $b=y$ or $b=xy$ (when $k$ is even).\\
\end{rmq}

We now prove a more general version of Theorem \hyperlink{ThD}{D} : we are able to observe monotonicity in the evolution of the bias in the subextensions of $\mathcal{D}_n/\Q$ and $\mathcal Q_n^+/\Q$.

\begin{thm}\label{Mono} Assume $\GRH$ and $\LIP$ ($\LI$ is sufficient for the dihedral case). For any $n \geq 3$ and $3 \leq i \leq n$, let $\mathcal{D}_n^{(i)} = \mathcal{D}_n^{\langle r_i, s\rangle}$ as in section 4.2, and $(\mathcal Q_n^+)^{(i)} = (\mathcal Q_n^+)^{\langle x_i, y\rangle}$ as in section 4.1. Then for any $\varepsilon > 0$ and any sufficiently large $n$, for $3 \leq i < j \leq n$ such that $i \leq n\frac{1 + \varepsilon}{2}$ and $j \geq n\left(\frac{1 + 3\varepsilon}{2}\right)$, we have $$\delta(\mathcal{D}_n/\mathcal{D}_n^{(j)}, C_1^{(j)}, C_{-1}^{(j)}) < \delta(\mathcal{D}_n/\mathcal{D}_n^{(i)}, C_1^{(i)}, C_{-1}^{(i)}),$$ $$1-\delta(\mathcal Q_n^+/(\mathcal Q_n^+)^{(j)}, C_1^{(j)}, C_{-1}^{(j)}) < 1- \delta(\mathcal Q_n^+/(\mathcal Q_n^+)^{(i)}, C_1^{(i)}, C_{-1}^{(i)})$$ and $$\delta(\mathcal Q_n^-/(\mathcal Q_n^-)^{(j)}, C_1^{(j)}, C_{-1}^{(j)}) < \delta(\mathcal Q_n^-/(\mathcal Q_n^-)^{(i)}, C_1^{(i)}, C_{-1}^{(i)}).$$
\end{thm}

\begin{demo} We combine the bounds of Propositions \ref{VarD}, \ref{EspD} and \ref{ConsD}, with Theorems \ref{DevFJ} and \ref{Dev}. With notations from Theorem \ref{Dev}, it is easy to see from Corollary \ref{Den} that $R = \{\chi_2, \chi_3\}$ so that $b_1=b_2=2, M=1$, and $\frac{b_3(C_1, C_{-1})}{b_4(C_1, C_{-1})}=1$. We thus find $$c_1\exp(-c_22^{2i-n}) < \delta(\mathcal{D}_n/\mathcal{D}_n^{(i)}, C_1^{(i)}, C_{-1}^{(i)}) < \exp\left(-c_3 \frac{2^{2i-n}}{n}\right)$$ for any $3 \leq i \leq n$ and for some absolute $c_1, c_2, c_3 > 0$.

In order to have $\delta(\mathcal{D}_n/\mathcal{D}_n^{(j)}, C_1^{(j)}, C_{-1}^{(j)}) < \delta(\mathcal{D}_n/\mathcal{D}_n^{(i)}, C_1^{(i)}, C_{-1}^{(i)}),$ it is therefore enough to have $$\exp\left(-c_3 \frac{2^{2j-n}}{n}\right) < c_1\exp(-c_22^{2i-n}).$$ This is equivalent to $$c_22^{2i} < c_3 \frac{2^{2j}}{n} + 2^n\log(c_1).$$ Now if $j \geq n\left(\frac{1 + 3\varepsilon}{2}\right)$ and $n$ is large enough then $c_3 \frac{2^{2j}}{n} \geq c_3 \frac{2^{n(1+3\varepsilon)}}{n} > 2^{n (1+2\varepsilon)}$ while if $i \leq n\left(\frac{1}{2} + \varepsilon\right)$ and $n$ is large enough we have $c_22^{2i} + 2^n \log(1/c_1) \leq 2^{n(1+2\varepsilon)}$ and the desired inequality holds.

The proof is similar in the case of $\mathcal Q_n^{\pm}$ by using the bounds of Proposition \ref{VarQ}, Proposition \ref{EspQ} and Proposition \ref{ConsQ} with Theorems \ref{DevFJ} and \ref{Dev}.
\end{demo}

\begin{rmq} A similar proof can be applied to the other extremely biased Chebotarev races from Theorems \ref{TabD} and \ref{TabQ}, which yields similar monotonicity results.
\end{rmq}

\subsection*{Acknowledgements}

The author would like to thank Florent Jouve for his doctoral advising and for suggesting to study quaternion extensions of $\Q$ in the context of Chebyshev's bias, Daniel Fiorilli for fruitful discussions on the subject, Philippe Cassou-Noguès for introducing him to the paper \cite{Fro1}, which was the starting point of this work, Gerhard Niklasch for technical help in bounding discriminants in the proof of Proposition \ref{ConsD}, Alain Debreil for providing \LaTeX \,code for the lattice of subgroups in Section 3, and the two anonymous referees for their suggestions and remarks which helped in improving the text.

\bibliographystyle{alpha}
\bibliography{bibli}

\begin{bibdiv}
\begin{biblist}

\bib{Arm}{article}{
      author={Armitage, J.~V.},
       title={Zeta functions with a zero at {$s={1\over 2}$}},
        date={1972},
     journal={Invent. Math.},
      volume={15},
       pages={199\ndash 205},
}

\bib{Cha}{article}{
      author={Cha, B.},
       title={Chebyshev's bias in function fields},
        date={2008},
     journal={Compos. Math.},
      volume={144},
      number={6},
       pages={1351\ndash 1374},
}

\bib{ChaIm}{article}{
      author={Cha, B.},
      author={Im, B.-H.},
       title={Chebyshev's bias in {G}alois extensions of global function
  fields},
        date={2011},
     journal={J. Number Theory},
      volume={131},
      number={10},
       pages={1875\ndash 1886},
}

\bib{DaMa}{article}{
      author={Damey, P.},
      author={Martinet, J.},
       title={Plongement d'une extension quadratique dans une extension
  quaternionienne},
        date={1973},
     journal={J. Reine Angew. Math.},
      volume={262/263},
       pages={323\ndash 338},
        note={Collection of articles dedicated to Helmut Hasse on his
  seventy-fifth birthday},
}

\bib{Deb}{book}{
      author={Debreil, A.},
       title={Groupes finis et treillis de leurs sous-groupes},
      series={Math{\'e}matiques en devenir},
   publisher={Calvage et Mounet},
        date={2016},
}

\bib{Dev}{inproceedings}{
      author={Devin, L.},
       title={Chebyshev’s bias for analytic {L}-functions},
organization={Cambridge University Press},
        date={2019},
   booktitle={Mathematical proceedings of the cambridge philosophical society},
       pages={1\ndash 38},
}

\bib{DevMen}{misc}{
      author={Devin, L.},
      author={Meng, X.},
       title={Chebyshev's bias for products of irreducible polynomials},
         how={\href{https://arxiv.org/pdf/1809.09662.pdf}{arXiv:1809.09662}},
        date={2020},
}

\bib{FiJo}{article}{
      author={Fiorilli, D.},
      author={Jouve, F.},
       title={Distribution of {F}robenius elements in families of {G}alois
  extensions},
        date={2020-01-15},
      eprint={http://arxiv.org/abs/2001.05428v1},
         url={http://arxiv.org/abs/2001.05428v1},
}

\bib{FiMa}{article}{
      author={Fiorilli, D.},
      author={Martin, G.},
       title={Inequities in the {S}hanks-{R}\'{e}nyi prime number race: an
  asymptotic formula for the densities},
        date={2013},
     journal={J. Reine Angew. Math.},
      volume={676},
       pages={121\ndash 212},
}

\bib{Lam4}{article}{
      author={Ford, K.},
      author={Harper, A.~J.},
      author={Lamzouri, Y.},
       title={Extreme biases in prime number races with many contestants},
        date={2019},
     journal={Math. Ann.},
      volume={374},
      number={1-2},
       pages={517\ndash 551},
}

\bib{Fro1}{article}{
      author={Fr\"{o}hlich, A.},
       title={Artin root numbers and normal integral bases for quaternion
  fields},
        date={1972},
     journal={Invent. Math.},
      volume={17},
       pages={143\ndash 166},
}

\bib{Fro2}{article}{
      author={Fr\"{o}hlich, A.},
       title={The {G}alois module structure of algebraic integer rings in
  fields with generalised quaternion group},
        date={1974},
       pages={81\ndash 86. Bull. Soc. Math. France M\'{e}m. 27},
}

\bib{Fro3}{book}{
      author={Fr\"{o}hlich, A.},
       title={Galois module structure of algebraic integers},
      series={Ergebnisse der Mathematik und ihrer Grenzgebiete (3) [Results in
  Mathematics and Related Areas (3)]},
   publisher={Springer-Verlag, Berlin},
        date={1983},
      volume={1},
}

\bib{Lam3}{article}{
      author={Harper, A.~J.},
      author={Lamzouri, Y.},
       title={Orderings of weakly correlated random variables, and prime number
  races with many contestants},
        date={2018},
     journal={Probab. Theory Related Fields},
      volume={170},
      number={3-4},
       pages={961\ndash 1010},
}

\bib{Isa}{book}{
      author={Isaacs, I.~M.},
       title={Character theory of finite groups},
   publisher={Dover Publications, Inc., New York},
        date={1994},
        note={Corrected reprint of the 1976 original [Academic Press, New York;
  MR0460423 (57 \#417)]},
}

\bib{IwKo}{book}{
      author={Iwaniec, H.},
      author={Kowalski, E.},
       title={Analytic number theory},
      series={American Mathematical Society Colloquium Publications},
   publisher={American Mathematical Society, Providence, RI},
        date={2004},
      volume={53},
}

\bib{JeLeYu}{book}{
      author={Jensen, C.~U.},
      author={Ledet, A.},
      author={Yui, N.},
       title={Generic polynomials},
      series={Mathematical Sciences Research Institute Publications},
   publisher={Cambridge University Press, Cambridge},
        date={2002},
      volume={45},
        note={Constructive aspects of the inverse Galois problem},
}

\bib{Kac}{article}{
      author={Kaczorowski, J.},
       title={A contribution to the {S}hanks-{R}\'{e}nyi race problem},
        date={1993},
     journal={Quart. J. Math. Oxford Ser. (2)},
      volume={44},
      number={176},
       pages={451\ndash 458},
}

\bib{Kac2}{article}{
      author={Kaczorowski, J.},
       title={On the distribution of primes (mod {$4$})},
        date={1995},
     journal={Analysis},
      volume={15},
      number={2},
       pages={159\ndash 171},
}

\bib{KT1}{article}{
      author={Knapowski, S.},
      author={Tur\'{a}n, P.},
       title={Comparative prime-number theory. i-viii},
        date={1962-1963},
     journal={Acta Math. Acad. Sci. Hungar.},
      volume={13-14},
}

\bib{KT2}{article}{
      author={Knapowski, S.},
      author={Tur\'{a}n, P.},
       title={Further developments in the comparative prime-number theory.
  i-vi},
        date={1964-1966},
     journal={Acta Arith.},
      volume={9-12},
}

\bib{LMO}{article}{
      author={Lagarias, J.~C.},
      author={Montgomery, H.~L.},
      author={Odlyzko, A.~M.},
       title={A bound for the least prime ideal in the {C}hebotarev density
  theorem},
        date={1979},
     journal={Invent. Math.},
      volume={54},
      number={3},
       pages={271\ndash 296},
}

\bib{Lam2}{article}{
      author={Lamzouri, Y.},
       title={The {S}hanks-{R}\'{e}nyi prime number race with many
  contestants},
        date={2012},
     journal={Math. Res. Lett.},
      volume={19},
      number={3},
       pages={649\ndash 666},
}

\bib{Lam1}{article}{
      author={Lamzouri, Y.},
       title={Prime number races with three or more competitors},
        date={2013},
     journal={Math. Ann.},
      volume={356},
      number={3},
       pages={1117\ndash 1162},
}

\bib{Lan}{book}{
      author={Lang, S.},
       title={Algebraic number theory},
     edition={Second},
      series={Graduate Texts in Mathematics},
   publisher={Springer-Verlag, New York},
        date={1994},
      volume={110},
}

\bib{MN}{article}{
      author={Martin, G.},
      author={Ng, N.},
       title={Inclusive prime number races},
        date={2020},
     journal={Trans. Amer. Math. Soc.},
      volume={373},
      number={5},
       pages={3561\ndash 3607},
}

\bib{MoOd}{article}{
      author={Montgomery, H.~L.},
      author={Odlyzko, A.~M.},
       title={Large deviations of sums of independent random variables},
        date={1988},
     journal={Acta Arith.},
      volume={49},
      number={4},
       pages={427\ndash 434},
}

\bib{MuMu}{book}{
      author={Murty, M.~Ram},
      author={Murty, V.~Kumar},
       title={Non-vanishing of {$L$}-functions and applications},
      series={Modern Birkh\"{a}user Classics},
   publisher={Birkh\"{a}user/Springer Basel AG, Basel},
        date={1997},
        note={[2011 reprint of the 1997 original] [MR1482805]},
}

\bib{Neu}{book}{
      author={Neukirch, J.},
       title={Algebraic number theory},
      series={Grundlehren der Mathematischen Wissenschaften [Fundamental
  Principles of Mathematical Sciences]},
   publisher={Springer-Verlag, Berlin},
        date={1999},
      volume={322},
        note={Translated from the 1992 German original and with a note by
  Norbert Schappacher, With a foreword by G. Harder},
}

\bib{Ng}{thesis}{
      author={Ng, N.},
       title={Limiting distributions and zeros of artin l-functions},
        type={Ph.D. Thesis},
        date={2000},
}

\bib{Od}{article}{
      author={Odlyzko, A.~M.},
       title={Bounds for discriminants and related estimates for class numbers,
  regulators and zeros of zeta functions: a survey of recent results},
        date={1990},
     journal={S\'{e}m. Th\'{e}or. Nombres Bordeaux (2)},
      volume={2},
      number={1},
       pages={119\ndash 141},
}

\bib{Pla}{article}{
      author={Plans, B.},
       title={On the minimal number of ramified primes in some solvable
  extensions of {$\Bbb Q$}},
        date={2004},
     journal={Pacific J. Math.},
      volume={215},
      number={2},
       pages={381\ndash 391},
}

\bib{RuSa}{article}{
      author={Rubinstein, M.},
      author={Sarnak, P.},
       title={Chebyshev's bias},
        date={1994},
     journal={Experiment. Math.},
      volume={3},
      number={3},
       pages={173\ndash 197},
}

\bib{RudSa}{article}{
      author={Rudnick, Z.},
      author={Sarnak, P.},
       title={Zeros of principal {$L$}-functions and random matrix theory},
        date={1996},
     journal={Duke Math. J.},
      volume={81},
      number={2},
       pages={269\ndash 322},
        note={A celebration of John F. Nash, Jr.},
}

\bib{Ser}{book}{
      author={Serre, J.-P.},
       title={Corps locaux},
   publisher={Hermann},
        date={1980},
}

\bib{SerCheb}{article}{
      author={Serre, Jean-Pierre},
       title={Quelques applications du th\'{e}or\`eme de densit\'{e} de
  {C}hebotarev},
        date={1981},
     journal={Inst. Hautes \'{E}tudes Sci. Publ. Math.},
      number={54},
       pages={323\ndash 401},
}

\bib{Ska}{thesis}{
      author={Skabelund, D.},
       title={Character {T}ables of {M}etacyclic {G}roups},
        type={Master's Thesis},
        date={2013},
}

\bib{Cheb}{article}{
      author={Tchebychev, P.},
       title={Lettre de {M}. le {P}rofesseur {T}chébychev à {M}. {F}uss sur
  un nouveau théorème relatif aux nombres premiers contenus dans les formes
  {$4n + 1$} et {$4n + 3$}},
        date={1853},
     journal={Bull. Classe Phys. Acad. Imp. Sci. St. Petersburg},
}

\bib{Toy}{inproceedings}{
      author={T{\^o}yama, H.},
       title={A note on the different of the composed field},
organization={Department of Mathematics, Tokyo Institute of Technology},
        date={1955},
   booktitle={Kodai mathematical seminar reports},
      volume={7},
       pages={43\ndash 44},
}

\bib{Win}{article}{
      author={Wintner, A.},
       title={On the distribution function of the remainder term of the prime
  number theorem},
        date={1941},
     journal={Amer. J. Math.},
      volume={63},
       pages={233\ndash 248},
}

\end{biblist}
\end{bibdiv}
\nocite{*}

\bigskip
\footnotesize
\noindent Alexandre Bailleul, \textsc{Univ. Bordeaux, IMB, UMR 5251, F 33405, Talence, France}\\
\textit{Email address:} \textsf{alexandre.bailleul@math.u-bordeaux.fr}

\end{document}